\documentclass[11pt,a4paper]{amsart}
\usepackage[T1]{fontenc}
\usepackage{lmodern}
\usepackage{amsmath,amssymb,amsthm,mathtools}
\usepackage[margin=29mm,headheight=14pt]{geometry}
\usepackage{xurl}
\usepackage{tikz}
\usepackage{xcolor}
\usetikzlibrary{decorations.markings}
\definecolor{figureblue}{RGB}{31,75,105}
\usepackage[dvipdfmx,unicode,hidelinks,bookmarksnumbered=true]{hyperref}
\hypersetup{
  pdftitle={Bennequin-sharp links are strongly quasipositive},
  pdfsubject={Bennequin sharpness, convex surfaces, and strong quasipositivity},
  pdfauthor={Nobuo Iida and Tatsumasa Suzuki},
  pdfkeywords={Bennequin inequality, strong quasipositivity, convex surfaces, sutured Floer homology}
}
\numberwithin{equation}{section}
\newtheorem{theorem}{Theorem}[section]
\newtheorem{lemma}[theorem]{Lemma}
\newtheorem{proposition}[theorem]{Proposition}
\newtheorem{corollary}[theorem]{Corollary}
\newtheorem{claim}{Claim}[section]
\newtheorem{question}[theorem]{Question}
\newtheorem{maintheorem}{Theorem}

\theoremstyle{definition}
\newtheorem{definition}[theorem]{Definition}
\theoremstyle{remark}
\newtheorem{remark}[theorem]{Remark}
\newcommand{\F}{\mathbb F_2}
\newcommand{\R}{\mathbb R}
\newcommand{\std}{\mathrm{std}}
\newcommand{\Int}{\operatorname{int}}
\newcommand{\slbar}{\overline{\operatorname{sl}}}
\newcommand{\slink}{\operatorname{sl}}

\newcommand{\SFH}{\operatorname{SFH}}
\newcommand{\EH}{\operatorname{EH}}
\newcommand{\HFhat}{\widehat{\operatorname{HF}}}
\newcommand{\Lie}{\mathcal L}
\title[Bennequin sharpness and strong quasipositivity]{Bennequin-sharp links are strongly quasipositive}

\author{Nobuo Iida}
\address{Kavli Institute for the Physics and Mathematics of the Universe (WPI), The University of Tokyo, 5-1-5 Kashiwanoha, Kashiwa, Chiba 277-8583, Japan}
\email{iidanobuo1224@g.ecc.u-tokyo.ac.jp}

\author{Tatsumasa Suzuki}
\address{Graduate School of Mathematical Sciences, The University of Tokyo, 3-8-1 Komaba, Meguro-ku, Tokyo 153-8914, Japan}
\email{suzukit@ms.u-tokyo.ac.jp}

\date{\today}

\begin{document}
\begin{abstract}
We give a contact-geometric proof that an oriented link in the standard contact three-sphere is strongly quasipositive if and only if its maximal total self-linking number equals minus the maximal Euler characteristic of a Seifert surface.
More precisely, every Seifert surface with positively transverse boundary satisfying the Bennequin equality is smoothly ambiently isotopic to a quasipositive surface.
This answers Hedden's Problem 1.69 and Problem 1.70(b) in the K3 problem list.
\end{abstract}

\maketitle

\section{Introduction}
Let $L\subset S^3$ be a nonempty smooth oriented link.
A Seifert surface for $L$ is a compact oriented embedded surface $F \subset S^3$ with $\partial F=L$ and no closed components. We allow $F$ to be disconnected and set
\begin{equation}\label{eq:chi-max}
\chi_{\max}(L)=\max_{\partial F=L}\chi(F).
\end{equation}
To see that the maximum exists, let $l$ be the number of components of $L$, let $c$ be the number of components of $F$, and let $g$ be the sum of their genera.
Each component of $F$ has at least one boundary component, so $c\leq l$ and $\chi(F)=2c-2g-l \leq l$.
The possible Euler characteristics thus form a nonempty set of integers bounded above.
\par
For a positive transverse link $T=\partial F$ in a contact three-manifold, write $\slink_F(T)$ for its total self-linking number relative to $F$.
In the standard contact three-sphere $(S^3, \xi_{\std})$ this number is independent of the choice of Seifert surface, and we also write $\slink(T)$.
Let $\slbar(L)$ denote the maximum of $\slink(T)$ over positive transverse representatives of $L$ in $(S^3, \xi_{\std})$.
The Bennequin inequality is 
\begin{equation}\label{eq:bennequin}
\slbar(L)\leq -\chi_{\max}(L).
\end{equation}
We call $L$ \emph{Bennequin-sharp} if equality holds in \eqref{eq:bennequin}. 
Notice that for a knot $K \subset S^3$, $\chi_{\max}(K)=1-2g_3(K)$.
\par
An oriented link is strongly quasipositive if it is the closure of a braid that is a product of positive embedded band generators.
Its associated positive band surface is called a quasipositive Seifert surface; this surface is allowed to be disconnected. 
Any surface smoothly ambiently isotopic to such a band surface is also called a quasipositive surface. 

Strong quasipositivity implies equality in \eqref{eq:bennequin}. 
Indeed, let $T$ be the positive transverse closure of a strongly quasipositive braid on $n$ strands with $m$ positive embedded bands, and let $F$ be its band surface. 
The $n$ disks and $m$ bands give $\chi(F)=n-m$. 
Each positive band generator has exponent sum one, so the self-linking formula for closed braids gives $\slink(T)=m-n=-\chi(F)$. By \eqref{eq:bennequin} and the definition of $\chi_{\max}(L)$, 
\[
 -\chi(F)=\slink(T)\le\slbar(L)\le-\chi_{\max}(L)\le-\chi(F).
\]
Thus all inequalities are equalities.
\par
The converse asks whether every Bennequin-sharp link is strongly quasipositive.
See \cite[Conjecture~2.5]{Stoimenow}, for example.
The fibered case follows from the work of Hedden and
Etnyre--Van Horn-Morris \cite{Hedden,EVHM}.
Feller--Lewark--Lobb \cite[Corollary E]{FLL} prove the converse for knots
admitting a canonical Seifert surface of minimal genus, and
Stoimenow \cite[Corollary 4.31]{Stoimenow} proves it for non-split
links admitting a canonical Seifert surface of maximal Euler characteristic.
Here a canonical Seifert surface is one obtained from an oriented link
diagram by Seifert's algorithm.

\begin{maintheorem}\label{main:surface}
Let $F \subset (S^3, \xi_{\std})$ be a compact oriented (not necessarily connected) Seifert surface with no closed components and positively transverse boundary $T$.
If 
\begin{equation}\label{eq:surface-sharp}
\slink(T)=-\chi(F),
\end{equation}
then $F$ is smoothly ambiently isotopic to a quasipositive surface.
The isotopy is not required to preserve the initially specified transverse isotopy class of $T$.
\end{maintheorem}

\begin{maintheorem}\label{main:link}
For every nonempty oriented link $L \subset S^3$, 
\[
 \slbar(L)=-\chi_{\max}(L)
 \quad\Longleftrightarrow\quad
 L\text{ is strongly quasipositive}.
\]
Moreover, if these conditions hold, every Seifert surface $F$ for $L$ satisfying $\chi(F)=\chi_{\max}(L)$ is quasipositive up to smooth ambient isotopy. (Split links and disconnected Seifert surfaces are included.)
\end{maintheorem}
Theorem~\ref{main:link} answers Hedden's Problem 1.69 in the K3 problem list \cite{K3}.
Theorem~\ref{main:surface} also answers Remark (3) of that problem.
Theorem~\ref{main:link} follows from Theorem~\ref{main:surface}.
\par
The assertion about every surface realizing \eqref{eq:chi-max} also gives
an affirmative answer to the question of Baker and Motegi recorded in
\cite[Question 2.2]{Stoimenow}, concerning the quasipositivity of all
minimal-genus Seifert surfaces of a strongly quasipositive link. Our
formulation uses maximal Euler characteristic and allows disconnected
Seifert surfaces. The corresponding result for canonical Seifert surfaces
of non-split strongly quasipositive links is proved in
\cite[Corollary 4.33]{Stoimenow}.
\begin{corollary}[Knots]\label{cor:knots}
For every knot $K \subset S^3$,
\[
 \slbar(K)=2g_3(K)-1
 \quad\Longleftrightarrow\quad
 K\text{ is strongly quasipositive}.
\]
Every minimal-genus Seifert surface of such a knot is quasipositive up to smooth ambient isotopy.
\end{corollary}
Hedden's Problem 1.70(b) in the K3 problem list \cite{K3} asks whether strong quasipositivity is characterized by quasipositivity together with equality of the Seifert and smooth slice genera. 
For links, we use the Euler-characteristic formulation described in Remark (3) to that problem. 
Set
\[
\chi_{4, \max}(L)=\max_{\partial S=L} \chi(S),
\]
where $S\subset B^4$ is a compact oriented smoothly and properly embedded surface with no closed components;
disconnected surfaces are allowed.
Here a link is quasipositive if it is the closure of a braid that is a product of conjugates of positive Artin generators.
\begin{corollary}[Quasipositivity and the slice bound]\label{cor:qp-genus}
A nonempty oriented link $L \subset S^3$ is strongly quasipositive if and only if it is quasipositive and 
\[
\chi_{4, \max}(L)=\chi_{\max}(L).
\]
In particular, for every knot $K \subset S^3$, 
\[
 K\text{ is strongly quasipositive}
 \quad\Longleftrightarrow\quad
 \begin{gathered}
 K\text{ is quasipositive},\\
 g_4(K)=g_3(K).
 \end{gathered}
\]
\end{corollary}
\begin{proof}
A quasipositive link is the closure of a braid that is a product of $k$ positive bands on $n$ strands. 
By Rudolph's construction \cite[\S2]{RudolphSlice}, it bounds a quasipositive braid surface $S\subset B^4$ with $\chi(S)=n-k$. 
The transverse self-linking formula for closed braids then gives
\[
\slink(T)=k-n=-\chi(S),
\]
where $T$ is the positive transverse closure of the braid. 
The slice-Bennequin inequality gives
\[
-\chi(S)=\slink(T)\leq -\chi_{4, \max}(L)\leq -\chi(S),
\]
so $\slink(T)=-\chi_{4, \max}(L)$.
If $\chi_{4, \max}(L)=\chi_{\max}(L)$, Theorem~\ref{main:link} implies that $L$ is strongly quasipositive.
\par
Conversely, for a positive band surface $F$ of a strongly quasipositive link, its transverse boundary satisfies 
\[
 -\chi(F)=\slink_F(\partial F)
 \le-\chi_{4,\max}(L)\le-\chi_{\max}(L)\le-\chi(F).
\]
All inequalities are equalities. The knot statement follows from
$\chi_{4,\max}(K)=1-2g_4(K)$ and $\chi_{\max}(K)=1-2g_3(K)$.
\end{proof}
\subsection{Outline of the proof of Theorem~\ref{main:surface}}
We begin with a Seifert surface $F$ for a transverse link
$T$ in $(S^3,\xi_{\std})$ realizing the Bennequin equality.
The standard contact structure $\xi_{\std}$ induces a
characteristic foliation on $F$, whose singularities are
precisely the points where $\xi_{\std}$ is tangent to $F$.
Such a singularity is called positive or negative according to whether the orientation of $\xi_{\std}$ agrees or disagrees with that of $F$.
\par
The Bennequin equality and tightness allow us to cancel all negative singularities. To describe the resulting configuration, recall that a surface $F$ is convex if it admits a contact vector field $X$ transverse to $F$. Writing $\xi_{\std}=\ker \alpha$, its dividing set and positive and negative regions are
\[
\Gamma_F=\{p \in F \mid \alpha(X)(p)=0\}, \qquad R_\pm(F)=\{p \in F \mid \pm \alpha(X)(p)\geq 0\}. 
\]
After the cancellation and the convex-surface preparation in Section 2, the negative regions are annuli of two types: boundary collars of $T$ and interior annuli surrounding attracting closed leaves of the characteristic foliation.
Our main task is to eliminate the interior negative annuli without changing the smooth isotopy class of $F$.
\par
We achieve this by modifying the contact structure in a solid torus neighborhood of each selected interior annulus and replacing the annulus by a smoothly isotopic one.
The contact structure can be glued to the unchanged exterior.
\par
The two main ingredients are stated in Section 3, where we combine them to eliminate all interior negative annuli and deduce the main theorems.
\par
The first ingredient is a local replacement construction. It uses the six-point relation for the Honda--Kazez--Mati\'c
contact invariant of contact $3$-manifolds with convex
boundary. 
Three contact fillings are glued to the same exterior, and the contact invariant of the original filling is the sum of those of the two alternatives. 
Consequently, if the original contact invariant is nonzero, at least one alternative also has nonzero contact invariant.
Both alternatives admit suitable replacement annuli.
The construction of these fillings, including their common boundary marking, is given in Section 4.
\par
The second ingredient prepares a selected interior negative annulus for replacement. 
More precisely, we modify the characteristic one-form on $F$ to realize a prescribed rotation model near the annulus. The modified one-form must be extended over the adjacent positive regions without violating the contact condition. This is possible under the outlet condition, which ensures the positivity of the boundary integrals required by the area-form extension lemma.

If the outlet condition fails, we first perform a fold along a nonseparating curve in the obstructing positive region. 
This creates an auxiliary negative annulus and provides the missing outlets for the original annulus.
After replacing the original annulus, the auxiliary annulus also satisfies the outlet condition and can be replaced in turn. 
Thus two replacements eliminate both the original and auxiliary annuli. 

Section~5 proves the realization statement under the outlet condition, while Section~3 establishes the fold-and-replacement argument.
\par
Repeating these operations removes all interior negative annuli. 
The key step in passing from this configuration to a Legendrian ribbon is Hayden's ribbon-flexibility lemma \cite[Lemma~3.1]{Hayden}. 
Indeed, as proved in Lemma~\ref{lem:terminal}, the characteristic foliation can then be made Morse--Smale, with no negative singularities or attracting closed leaves. 
Hayden's lemma therefore gives an isotopy to a Legendrian ribbon in the final contact structure. 

Since the contact invariant remains nonzero, this contact structure is tight. By the uniqueness of the tight contact structure on $S^3$ and Gray stability, the ribbon can be carried to $(S^3,\xi_{\std})$. 
The correspondence between Legendrian ribbons and quasipositive surfaces established by Baader--Ishikawa and Hayden \cite{BI,Hayden} then shows that the original surface is smoothly isotopic to a quasipositive surface. 

\subsection*{Acknowledgements}
NI was supported by JSPS KAKENHI Grant Number JP26K16990 (Grant-in-Aid for Early-Career Scientists) and by the Start-up Fund of the Kavli Institute for the Physics and Mathematics of the Universe (Kavli IPMU), The University of Tokyo. 
TS was supported by JSPS KAKENHI, Grant Numbers JP26KJ0062 (Grant-in-Aid for JSPS Fellows) and JP26K16991 (Grant-in-Aid for Early-Career Scientists). 

\subsection*{AI Declaration}

\par
ChatGPT proposed a proof of  Proposition~\ref{prop:normal-form}, Proposition~\ref{prop:P}, Proposition~\ref{prop:L},  and Theorem~\ref{thm:elimination}.
The authors also used ChatGPT to assist with exploratory bibliographic searches, English-language editing, and figure generation. 
The authors independently checked all references and verified all mathematical arguments and computations, and take full responsibility for the content of the manuscript. 

\section{Convex surfaces and the annular normal form}\label{sec:setup}
Let $F$ be a compact oriented surface embedded in a positive cooriented contact three-manifold $(M, \xi=\ker \alpha)$. 
We give $\partial F$ its boundary orientation. 
The boundary $\partial F$ is said to be positively transverse if $\alpha$ is positive on positively oriented tangent vectors to $\partial F$. 
A vector field $X$ is called a contact vector field if 
\[
\mathcal{L}_X\alpha=h\alpha
\]
for some smooth function $h$, i.e. its local flow preserves $\xi$. 
The surface $F$ is convex if such a field $X$ is defined near $F$ and is transverse to it. 
Choose $X$ positively transverse to $F$ and put $H=\alpha(X)|_F$. 
We use the closed positive and negative regions
\[
\Gamma_F=H^{-1}(0), \qquad R_+(F)=\{H \geq 0\}, \quad R_-(F)=\{H \leq 0\}.
\]

\subsection{Disconnected surfaces and total self-linking}
Since we allow disconnected surfaces, we first observe that the Bennequin equality holds on each connected component whenever it holds for the entire surface. 

\begin{lemma}[Componentwise sharpness]\label{lem:component-sharpness}
Let $F=\bigsqcup_{a=1}^c F_a$ be an oriented surface with no closed components in a tight contact three-manifold, with positively transverse boundary $T$. 
Write $T_a=\partial F_a$. 
Then
\[
\slink_F(T)=\sum_{a=1}^c\slink_{F_a}(T_a).
\]
If $\slink_F(T)=-\chi(F)$, then
\begin{equation}\label{eq:component-sharpness}
 \slink_{F_a}(T_a)=-\chi(F_a)\qquad(1\le a\le c).
\end{equation}
\end{lemma}
\begin{proof}
Relative Euler numbers are additive over the disjoint components of the surface. 
The Bennequin inequality for each $F_a$ gives 
\[
d_a:=-\chi(F_a)-\slink_{F_a}(T_a)\geq 0.
\]
The sharpness assumption and the additivity of Euler characteristic imply $\sum_a d_a=0$. Hence each $d_a$ vanishes.
\end{proof}

\subsection{Protected annular configurations}

We next introduce the annular configuration used in our elimination procedure. As we will show below, a surface satisfying the Bennequin equality can be brought into this form by cancelling the negative singularities of its characteristic foliation and arranging the remaining negative regions as annuli. We also specify the boundary collars whose characteristic foliations will be kept fixed throughout the construction.
\par

In the definition below, we orient the characteristic foliation on an oriented
surface $F$ by the vector field $W$ determined by
\[
  \iota_W\Omega=\alpha|_F,
\]
where $\alpha$ is a positive contact form for $\xi$ and $\Omega$ is a positive area form on $F$.
This orientation is independent of the choices of $\alpha$ and $\Omega$.

\begin{definition}\label{def:protected}(Protected annular configuration).
Let $(M, \xi)$ be a positive cooriented contact three-manifold, and let $F\subset M$ be a compact oriented convex surface, possibly disconnected, with no closed components and with positively transverse boundary $T=T_1\sqcup\cdots \sqcup T_\ell$.
Choose a contact vector field on a neighborhood of $F$ that is positively transverse to $F$. 
A \emph{protected annular configuration} on $F$ consists of a closed annular collar $N_i$ of each $T_i$ and closed annuli $A_1, \dots, A_m \subset \operatorname{int}F$. 
The $N_i$ are pairwise disjoint, and $A_j$ are pairwise disjoint and disjoint from every $N_i$. 
Write $\partial N_i=T_i \sqcup C_i$, $N=\bigsqcup_i N_i$, and $C=\bigsqcup_i C_i$. 
We require
\begin{equation}
R_-(F)=N \sqcup A_1 \sqcup \cdots \sqcup A_m, \quad \Gamma_F=C \sqcup \partial A_1 \sqcup \cdots \sqcup \partial A_m.
\end{equation}\label{eq:normal-form}
On each $N_i$ the characteristic foliation is a product by arcs oriented from $C_i$ to $T_i$. 
These $N_i$'s are called the protected boundary collars.
\end{definition}

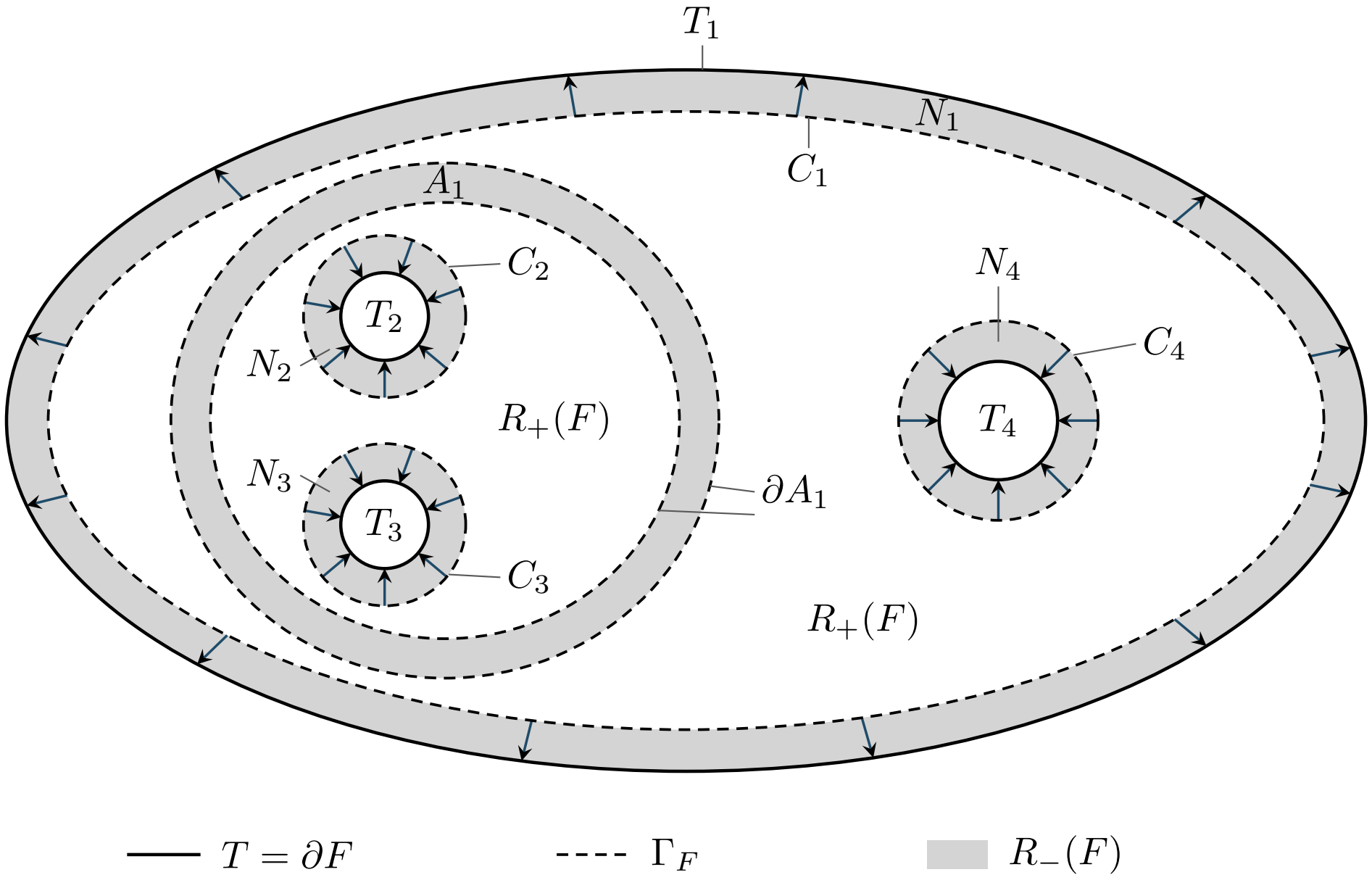
\begin{figure}[tbp]
\centering
\begingroup
\definecolor{pacflow}{RGB}{31,75,105}
\begin{tikzpicture}[
  x=1cm,y=1cm,font=\small,
  boundary/.style={draw=black,line width=.85pt},
  dividing/.style={draw=black,line width=.70pt,
    dash pattern=on 3pt off 2pt},
  leaf/.style={draw=pacflow,line width=.65pt,->,>=stealth},
  guide/.style={draw=black!65,line width=.40pt},
  every node/.style={inner sep=1.5pt}
]
\fill[black!17,even odd rule]
  (0,0) ellipse [x radius=6.20,y radius=3.20]
  (0,0) ellipse [x radius=5.82,y radius=2.82];

\fill[black!17,even odd rule]
  (-2.20,0) ellipse [x radius=2.50,y radius=2.35]
  (-2.20,0) ellipse [x radius=2.14,y radius=1.99];

\foreach \y in {-.95,.95}{
  \fill[black!17,even odd rule]
    (-2.75,\y) circle [radius=.74]
    (-2.75,\y) circle [radius=.40];
}
\fill[black!17,even odd rule]
  (2.85,0) circle [radius=.91]
  (2.85,0) circle [radius=.54];

\foreach \a in {12,40,80,100,134,166,194,224,256,286,320,348}{
  \draw[leaf]
    ({5.82*cos(\a)},{2.82*sin(\a)}) --
    ({6.20*cos(\a)},{3.20*sin(\a)});
}
\foreach \y in {-.95,.95}{
  \foreach \a in {20,70,120,170,220,270,320}{
    \draw[leaf]
      ({-2.75+.74*cos(\a)},{\y+.74*sin(\a)}) --
      ({-2.75+.40*cos(\a)},{\y+.40*sin(\a)});
  }
}
\foreach \a in {0,45,135,180,225,270,315}{
  \draw[leaf]
    ({2.85+.91*cos(\a)},{.91*sin(\a)}) --
    ({2.85+.54*cos(\a)},{.54*sin(\a)});
}

\draw[dividing] (0,0) ellipse [x radius=5.82,y radius=2.82];
\draw[dividing] (-2.20,0) ellipse [x radius=2.50,y radius=2.35];
\draw[dividing] (-2.20,0) ellipse [x radius=2.14,y radius=1.99];
\foreach \y in {-.95,.95}
  \draw[dividing] (-2.75,\y) circle [radius=.74];
\draw[dividing] (2.85,0) circle [radius=.91];

\draw[boundary] (0,0) ellipse [x radius=6.20,y radius=3.20];
\foreach \y in {-.95,.95}
  \draw[boundary] (-2.75,\y) circle [radius=.40];
\draw[boundary] (2.85,0) circle [radius=.54];

\node (t1) at (.15,3.63) {$T_1$};
\draw[guide] (t1.south) -- (.15,3.20);
\node at (2.30,2.79) {$N_1$};
\node (c1) at (1.12,2.28) {$C_1$};
\draw[guide] (c1.north) -- (1.12,2.767);

\node at (-2.75,.95) {$T_2$};
\node at (-2.75,-.95) {$T_3$};
\node at (2.85,0) {$T_4$};
\node (n2) at (-3.80,.50) {$N_2$};
\draw[guide] (n2.east) -- (-3.25,.65);
\node (n3) at (-3.80,-.50) {$N_3$};
\draw[guide] (n3.east) -- (-3.25,-.65);
\node (c2) at (-1.43,1.43) {$C_2$};
\draw[guide] (c2.west) -- (-2.166,1.404);
\node (c3) at (-1.43,-1.43) {$C_3$};
\draw[guide] (c3.west) -- (-2.166,-1.404);
\node (n4) at (2.85,1.43) {$N_4$};
\draw[guide] (n4.south) -- (2.85,.72);
\node (c4) at (4.37,.70) {$C_4$};
\draw[guide] (c4.west) -- (3.535,.599);

\node at (-2.20,2.17) {$A_1$};
\node at (-1.20,0) {$R_+(F)$};
\node at (1.62,-1.83) {$R_+(F)$};
\node (da) at (1.00,-.65) {$\partial A_1$};
\draw[guide] (da.west) -- (.218,-.597);
\draw[guide] (da.south west) -- (-.250,-.820);

\begin{scope}[yshift=-3.96cm]
  \draw[boundary] (-5.10,0) -- (-4.43,0);
  \node[anchor=west] at (-4.30,0) {$T=\partial F$};
  \draw[dividing] (-1.18,0) -- (-.51,0);
  \node[anchor=west] at (-.38,0) {$\Gamma_F$};
  \fill[black!17] (2.20,-.13) rectangle (2.75,.13);
  \node[anchor=west] at (2.88,0) {$R_-(F)$};
\end{scope}
\end{tikzpicture}
\endgroup
\caption{An example of a protected annular configuration on a
four-holed sphere. The shaded annuli are the protected boundary collars
$N_1,\ldots,N_4$ and the interior annulus $A_1$, whose union is
$R_-(F)$. The dashed curves form
$\Gamma_F=C_1\sqcup C_2\sqcup C_3\sqcup C_4\sqcup\partial A_1$;
the solid curves are the actual boundary components $T_i$.
The interiors of the positive regions are unshaded, and the three
inner disks are omitted from $F$. Arrows indicate the characteristic
leaves on the protected collars, oriented from $C_i$ to $T_i$.}
\label{fig:protected-configuration}
\end{figure}

The isotopies used in the annulus-elimination procedure are fixed near $T$. However, the final ambient isotopy in Theorem~\ref{main:surface} need not fix $T$, since the Gray isotopy identifying the resulting tight contact structure with $\xi_{\std}$ may move the boundary. 

\par
Tightness and Bennequin sharpness enter in the following existence statement.
\begin{proposition}[Protected annular form]\label{prop:normal-form}
Let $F$ be a compact oriented surface with no closed components in a tight contact three-manifold. 
Suppose its boundary $T$ is positively transverse and $\slink_F(T)=-\chi(F)$. 
Then an isotopy fixed near $\partial F$ puts $F$ into protected annular form. 
\end{proposition}

\begin{proof}
By Lemma~\ref{lem:component-sharpness}, it suffices to consider each connected component of $F$ separately. 
We first arrange the characteristic foliation and then construct a transverse contact vector field with the required dividing set. 

Let $\alpha$ be a positive contact form for $\xi$, and let $\Omega$ be a positive area form on $F$.
We orient the characteristic foliation by the vector field $W$ determined by
\[
  \iota_W\Omega=\alpha|_{TF}. 
\]
Since $\partial F$ is positively transverse, $W$ points strictly outwards along $\partial F$. Indeed, if $\tau$ is a positively oriented tangent vector to $\partial F$ and $\nu$ is an outward-pointing normal vector in $TF$, then
\[
  \Omega(\nu,\tau)>0,\qquad
  \Omega(W,\tau)=\alpha(\tau)>0. 
\]
It follows that the outward component of $W$ is positive. 

We therefore choose disjoint, sufficiently small one-sided collars of the boundary components on which the characteristic foliation consists of arcs directed towards $\partial F$. 
We keep these collars fixed in the following isotopies. 

By the relative genericity theorem for characteristic foliations, we may perturb the interior of $F$ so that its characteristic foliation is Morse--Smale, while preserving the chosen boundary collars; compare \cite[\S3]{EVHM}
\footnote{%
A point $p$ of a flow $\varphi_t$ is \emph{nonwandering} if, for every neighborhood $U$ of $p$ and every $\tau>0$, there exist $q\in U$ and $t>\tau$ such that $\varphi_t(q)\in U$, with the trajectory defined throughout $[0,t]$. 
The \emph{nonwandering set} consists of all such points. 
A singularity of a vector field $W$ is \emph{hyperbolic in the dynamical sense} if every eigenvalue of its
linearization has nonzero real part. 
A periodic orbit is hyperbolic if the derivative of its Poincar\'e first-return map has modulus different from one. 
A flow is \emph{Morse--Smale} if its nonwandering set consists of finitely many hyperbolic singularities
and periodic orbits, and their stable and unstable manifolds intersect transversely. 
On a surface, this excludes saddle connections. 
This dynamical use of ``hyperbolic'' differs from the elliptic/hyperbolic classification of characteristic
singularities. 
In the latter terminology, elliptic singularities have index $+1$, while hyperbolic singularities are saddles of index $-1$. 
In a Morse--Smale characteristic foliation, an elliptic singularity is a source or sink and is also hyperbolic in the dynamical sense. 
}. 

Let $e_\pm$ and $h_\pm$ denote the numbers of
positive/negative elliptic and positive/negative
hyperbolic singularities, respectively. The
Poincar\'e--Hopf formula and the self-linking formula
give
\[
  \chi(F)=(e_++e_-)-(h_++h_-),
\]
and
\[
  \slink_F(T)=-(e_+-h_+)+(e_--h_-);
\]
see \cite[Lemma~2.3 and its proof]{EVHM}.
Consequently,
\[
  -\chi(F)-\slink_F(T)=2(h_--e_-).
\]
The Bennequin equality therefore implies
\begin{equation}\label{eq:negative-count}
  e_-=h_-.
\end{equation}

We next eliminate the negative singularities.
We use the cancellation argument in
\cite[Lemma~3.3]{EVHM}.
First, any repelling periodic orbit can be broken by
introducing a pair consisting of a positive elliptic
and a positive hyperbolic singularity. This does not
change either $e_-$ or $h_-$.

For each negative elliptic singularity, consider its
basin of attraction. The argument in the cited lemma
uses tightness to show that the boundary of this
basin must contain a negative hyperbolic singularity
with which the elliptic singularity can be cancelled.
Otherwise, the elimination of the positive
singularities along the basin would produce an
overtwisted disk. Repeating the cancellation removes
all negative elliptic singularities.

Each cancellation removes one negative elliptic and
one negative hyperbolic singularity, so
\eqref{eq:negative-count} remains valid. Once
$e_-=0$, it follows that $h_-=0$ as well.
Thus no negative singularities remain. All these
operations are supported away from $\partial F$.
The cancellation argument uses only tightness and
the Bennequin equality, not the fibering assumption
in the ambient theorem of \cite{EVHM}.

We retain the Morse--Smale condition after the
cancellations, making a further sufficiently small
interior perturbation if necessary. Such a perturbation
can be chosen without introducing negative
singularities and without changing the characteristic
foliation on the fixed boundary collars.

We now construct the dividing set, following
\cite[proof of Lemma~2.1]{EVHM}.
Since there are no negative singularities, choose
disjoint annular neighborhoods of the attracting
periodic orbits, together with the one-sided product
collars of $\partial F$ chosen above. Denote their
union by $U_-$.

Next, choose a neighborhood $U_+$ of the positive
elliptic singularities, the stable separatrices of
positive hyperbolic singularities, and the repelling
periodic orbits.

As explained in the proof of Lemma~2.1 of
\cite{EVHM}, these neighborhoods can be chosen
disjointly, with boundaries transverse to the
characteristic foliation, so that
\[
  F\setminus(U_+\cup U_-)
\]
is a disjoint union of product annuli foliated by
arcs running from $U_+$ to $U_-$.

Choose a core circle transverse to the characteristic
foliation in each product annulus, and let $\Gamma$
be the union of these circles. By the construction
in \cite[proof of Lemma~2.1]{EVHM}, the multicurve
$\Gamma$ divides the characteristic foliation in
the sense of Giroux.\footnote{%
A multicurve $\Gamma$ divides an oriented singular
foliation $\mathcal F$ on an oriented surface if
there exist a one-form $\beta$ representing
$\mathcal F$ and a smooth function $f$, with
$\Gamma=f^{-1}(0)$ a regular level set, such that
\[
  f\,d\beta+\beta\wedge df>0.
\]
Equivalently, $\ker(f\,dt+\beta)$ defines an
$I$-invariant contact structure on a neighborhood
of the surface, inducing $\mathcal F$ as its
characteristic foliation and $\Gamma$ as its
dividing set; see \cite[\S3.1.4]{Honda}.}
In particular, Lemma~2.1 of \cite{EVHM}
provides a contact vector field $X$ positively
transverse to $F$ whose dividing set is $\Gamma$.

We identify the positive and negative regions
determined by $X$. In an invariant flow neighborhood,
we may write
\[
  \alpha=f\,dt+\beta,\qquad X=\partial_t,
\]
where
\[
  f\,d\beta+\beta\wedge df>0,
  \qquad \Gamma=f^{-1}(0).
\]
Here the contact form has been positively rescaled
to be invariant under the flow of $X$.

Let $\Omega$ be a positive area form on $F$, and
define $W$ by $\iota_W\Omega=\beta$. Along $\Gamma$,
the contact condition gives
\[
  0<\beta\wedge df=-W(f)\Omega.
\]
Thus $W(f)<0$ on $\Gamma$. In other words, the
oriented characteristic foliation crosses $\Gamma$
from the positive region to the negative region.

Since the characteristic leaves in the product
annuli run from $U_+$ to $U_-$, the latter lies
on the negative side of the chosen dividing curves.
Consequently, the negative region is obtained by
enlarging the components of $U_-$ to the dividing
curves.

Each component of $U_-$ is either an annular
neighborhood of an attracting periodic orbit or
a one-sided product collar of a boundary component
$T_i$. Enlarging these neighborhoods through the
adjacent product annuli preserves their annular
topology. Hence every component of $R_-(F)$ is an
annulus.

Let $N_i$ denote the negative annulus containing
$T_i$, and let $C_i$ be its other boundary circle.
The characteristic foliation on $N_i$ consists of
arcs oriented from $C_i$ to $T_i$, since it extends
the original product foliation on the boundary
collar through a product annulus.

Denote the remaining components of $R_-(F)$ by
$A_1,\ldots,A_m$. These are annuli surrounding
attracting periodic orbits. We have therefore
obtained
\[
  R_-(F)=N\sqcup A_1\sqcup\cdots\sqcup A_m,
  \qquad N=\bigsqcup_iN_i,
\]
and
\[
  \Gamma_F=C\sqcup\partial A_1\sqcup\cdots
  \sqcup\partial A_m,
  \qquad C=\bigsqcup_iC_i.
\]
These are precisely the conditions in
\eqref{eq:normal-form}.

The isotopies used in the construction are fixed
near $\partial F$, and the characteristic foliations
on the resulting boundary collars are the prescribed
product foliations. For disconnected $F$, the same
construction is carried out independently on each
connected component.
\end{proof}

\begin{figure}[tbp]
\centering
\begin{tikzpicture}[
  x=1cm,y=1cm,font=\small,
  flow/.style={draw=figureblue,line width=.58pt},
  forward/.style args={#1}{postaction={decorate},decoration={markings,
    mark=at position #1 with {\arrow{stealth}}}},
  seam/.style={draw=black!55,densely dotted,line width=.5pt}]
  \fill[black!7] (0,-.82) rectangle (3.60,.82);
  \foreach \x in {0,3.60} \draw[seam] (\x,-.82) -- (\x,1.10);
  \foreach \x in {.40,1.10,1.80,2.50,3.20}
    \draw[flow,forward=.71] (\x,1.09) -- (\x,-.82);
  \draw[line width=.75pt] (0,.82) -- (3.60,.82)
    node[right=3pt] {$C_i$};
  \draw[line width=1.05pt] (0,-.82) -- (3.60,-.82)
    node[right=3pt] {$T_i$};
  \node at (1.80,-1.45) {Protected collar $N_i$};
  \begin{scope}[xshift=5.95cm]
    \fill[black!7] (0,-.82) rectangle (4.20,.82);
    \begin{scope}
      \clip (0,-1.10) rectangle (4.20,1.10);
      \foreach \side in {-1,1}
        \foreach \a/\pos in {1.50/.82,.60/.56,.24/.32}
          \draw[flow,forward=\pos] plot[domain=0:1,samples=32,variable=\s]
            ({4.20*\s},{\side*\a*exp(ln(.4)*\s)});
    \end{scope}
    \foreach \x in {0,4.20} \draw[seam] (\x,-1.10) -- (\x,1.10);
    \foreach \u in {-.82,.82}
      \draw[line width=.75pt] (0,\u) -- (4.20,\u)
        node[right=3pt] {$\partial A_j$};
    \draw[flow,line width=1.1pt,forward=.72] (0,0) -- (4.20,0)
      node[right=3pt,text=black] {$\ell$};
    \node at (2.10,-1.45) {Interior annulus $A_j$};
  \end{scope}
\end{tikzpicture}
\caption{The two kinds of negative annuli in the protected configuration.
In each panel the vertical sides are identified, the negative region is
shaded, and arrows orient the characteristic leaves. On $N_i$ the leaves
run from the dividing curve $C_i$ to the transverse boundary $T_i$.
On $A_j$ they spiral toward the closed leaf $\ell$ from both sides.}
\label{fig:protected-annuli}
\end{figure}
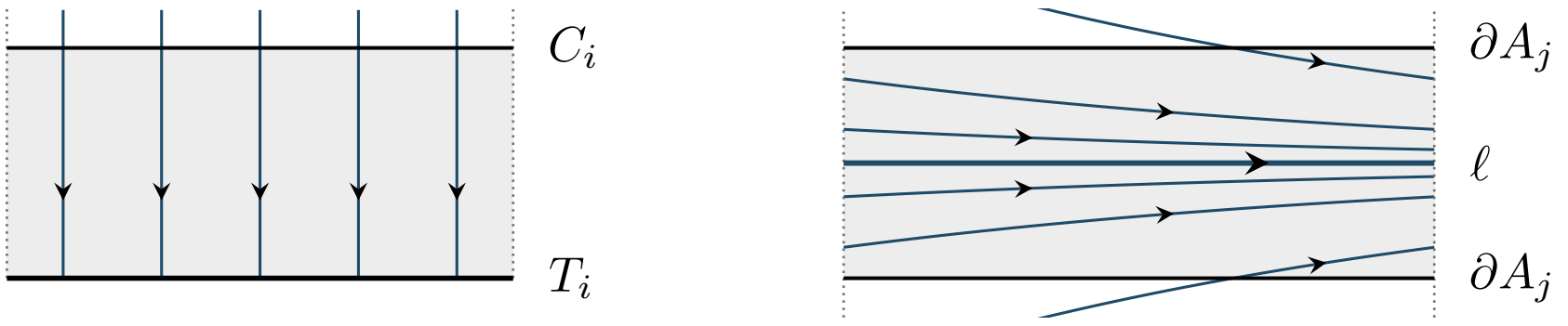

\begin{lemma}[Interior annular cores]
\label{lem:annular-core}
In a protected annular configuration in a tight
contact three-manifold, the core of an interior
negative annulus does not bound a disk in $F$.
\end{lemma}

\begin{proof}
Let $A\subset\operatorname{Int}F$ be an interior
negative annulus, and let $\gamma$ be a core circle
of $A$. By the definition of a protected annular
configuration,
\[
  \partial A=b_0\sqcup b_1\subset\Gamma_F.
\]
In particular,
\[
  \gamma\cap\Gamma_F=\varnothing.
\]

We first verify that $\gamma$ is nonisolating
with respect to $\Gamma_F$.

Identify $A$ with $S^1\times[0,1]$ so that
\[
  \gamma=S^1\times\{1/2\}.
\]
Cutting $A$ along $\gamma$ produces two annuli,
each having one boundary circle in $\Gamma_F$.
All other components of $F\setminus\Gamma_F$
are unchanged by this cutting. By the protected
annular configuration, each of these components
has a dividing curve in its boundary.
Consequently, the closure of every component of
\[
  F\setminus(\Gamma_F\cup\gamma)
\]
meets $\Gamma_F$. This is precisely the
nonisolating condition.

We next explain how to make $\gamma$ Legendrian
without changing the characteristic foliation
near $\partial F$.

The construction in the proof of Honda's
Legendrian realization principle
\cite[Theorem~3.7]{Honda} applies to a
nonisolating curve. It produces a singular
foliation $\mathcal F'$ divided by $\Gamma_F$
for which $\gamma$ is tangent to the foliation,
possibly passing through singular points.

In the present situation, $\gamma$ is contained
in the interior annulus $A$, and both components
of $\partial A$ belong to $\Gamma_F$. Hence
the construction can be performed while leaving
the original characteristic foliation unchanged
on every component of $F\setminus\Gamma_F$
meeting $\partial F$.

We may therefore apply the relative realization
theorem for convex surfaces with positively
transverse boundary
\cite[Theorem~2.2]{EVHM}.
It gives an isotopy of $F$, fixed near
$\partial F$, whose final characteristic
foliation is $\mathcal F'$.

Let $F'$ be the resulting convex surface and
let $\gamma'\subset F'$ be the image of
$\gamma$. Then $\gamma'$ is Legendrian.
Moreover, under the identification induced by
the isotopy, the dividing set is still
$\Gamma_F$. In particular,
\[
  \gamma'\cap\Gamma_{F'}=\varnothing.
\]

For a Legendrian curve on a convex surface,
the contact twisting relative to the surface
framing is
\[
  \operatorname{tw}(\gamma',F')
  =-\frac12\#(\gamma'\cap\Gamma_{F'});
\]
see \cite[\S3.3.1]{Honda}. Thus
\begin{equation}\label{eq:core-twisting}
  \operatorname{tw}(\gamma',F')=0.
\end{equation}

Suppose now, for a contradiction, that
$\gamma$ bounds an embedded disk $D\subset F$.
The isotopy carries $D$ to an embedded disk
$D'\subset F'$ with
\[
  \partial D'=\gamma'.
\]

Since $D'$ lies in $F'$, the surface framing
of $\gamma'$ induced by $F'$ agrees with
the disk framing induced by $D'$.
Equation~\eqref{eq:core-twisting} therefore
gives
\[
  \operatorname{tb}_{D'}(\gamma')=0.
\]

On the other hand, the Legendrian Bennequin
inequality in a tight contact three-manifold
implies
\[
  \operatorname{tb}_{D'}(\gamma')
  +|\operatorname{rot}_{D'}(\gamma')|
  \leq-\chi(D')=-1.
\]
In particular,
\[
  \operatorname{tb}_{D'}(\gamma')\leq-1,
\]
contradicting
$\operatorname{tb}_{D'}(\gamma')=0$.

Therefore $\gamma$ cannot bound an embedded
disk in $F$.
\end{proof}

\subsection{Legendrian ribbon}
\begin{lemma}[Relative genericity with positive divergence]
\label{lem:positive-genericity}
Let $Q$ be a compact oriented surface with boundary, $\omega$ a positive area form, and $Z$ a smooth vector field pointing strictly outwards on $\partial Q$. 
Suppose $\Lie_Z\omega=d(\iota_Z\omega)>0$. 
Then there is a smooth vector field $Z'$ arbitrarily close to $Z$ in $C^1$, equal to $Z$ near $\partial Q$, such that 
\[
 d(\iota_{Z'}\omega)>0
\]
and its interior dynamics are Morse--Smale: the nonwandering set consists of finitely many hyperbolic singularities and periodic orbits, and the stable and unstable manifolds intersect transversely. 
In particular, there are no saddle connections. 
Every periodic orbit of $Z'$ is repelling. 
\end{lemma}
\begin{proof}
Choose a boundary collar $C=\{0\le s<\varepsilon\}$, with $s$ increasing inwards, on which 
\[
Zs<-c<0. 
\]
Cap the boundary circles and extend $Z$ smoothly to the resulting closed oriented surface. 
By Peixoto's density theorem, structurally stable vector fields are dense in the $C^1$ topology on a compact surface. 
We may therefore approximate the extension arbitrarily closely in $C^1$ by a structurally stable vector field. 
Since structural stability is open in the $C^1$ topology and smooth vector fields are $C^1$-dense, the approximating vector field may be chosen smooth. 
By Peixoto's characterization of structural stability, the resulting vector field $\widetilde Z$ is Morse--Smale (see~\cite{Peixoto, Perko}). 

Fix a smooth cutoff $\chi(s)$ equal to $0$ near $s=0$ and to $1$ near $s=\epsilon$, and put
\[
Z'=(1-\chi)Z+\chi\widetilde Z
\quad\text{on }C,
\qquad
Z'=\widetilde Z\quad\text{off }C. 
\]
For a sufficiently close approximation, $Z's<-c/2$ on $C$. 
Thus $Z'$ has no singularities or periodic orbits in $C$, and every forward trajectory entering $C$ exits $Q$ through $\partial Q$. 
In particular, there are no nonwandering points in $C$. 

Any connecting orbit between interior critical elements cannot enter $C$, since every forward trajectory entering $C$ exits $Q$ through $\partial Q$ and cannot return.
Moreover, every periodic orbit is repelling, and every hyperbolic singularity is a source or a saddle because $\operatorname{div}_{\omega}Z'>0$.
Thus the only nonautomatic transversality condition concerns saddle-to-saddle connections. 
Such connections lie entirely in the region where $Z'=\widetilde Z$, so their transversality is inherited from the Morse--Smale field $\widetilde Z$. 

Since the cutoff is fixed, $Z'$ remains arbitrarily $C^1$-close to $Z$, so the strict inequality $d(\iota_{Z'}\omega)>0$ persists. 

For a periodic orbit $\gamma$ of period $T_\gamma$, the derivative of its transverse first-return map is 
\[
D\mathcal P_\gamma
=
\exp\left(
\int_0^{T_\gamma}
\operatorname{div}_{\omega}Z'
(\gamma(t))\,dt
\right)>1.
\]
Indeed, the derivative of the flow has multiplier $1$ in the orbit direction, and its area multiplier is the displayed exponential. 
Hence every periodic orbit of $Z'$ is repelling. 
\end{proof}

\begin{figure}[tbp]
\centering
\begin{tikzpicture}[
  x=1cm,y=1cm,font=\small,
  trajectory/.style={draw=black!60,line width=.48pt},
  core/.style={draw=figureblue,line width=1.15pt},
  oriented/.style={postaction={decorate},decoration={markings,
    mark=at position .72 with {\arrow{stealth}}}}]
   Both panels use polar coordinates, with central circle of radius 1.05.
  \foreach \offset in {0,5.75}{
    \begin{scope}[xshift=\offset cm]
      \fill[black!6,even odd rule] (0,0) circle (1.49) (0,0) circle (.61);
      \draw[line width=.55pt] (0,0) circle (1.49) (0,0) circle (.61);
    \end{scope}
  }
   A representative foliation near a repelling periodic orbit.
  \foreach \offset in {5,125,245}
    \foreach \side in {-1,1}
      \draw[trajectory,oriented] plot[domain=0:97,samples=32,variable=\t]
        ({\offset+\t}:{1.05+\side*.44*exp((\t-97)/54)});
  \draw[core,oriented] (0:1.05) arc (0:180:1.05);
  \draw[core,oriented] (180:1.05) arc (180:360:1.05);
  \node at (0,-1.82) {Repelling closed leaf};
  \draw[->,>=stealth,line width=.6pt] (1.90,0) -- (3.85,0);
   The two stable branches belong to one loop edge in this example.
  \begin{scope}[xshift=5.75cm]
    \foreach \side in {-1,1}
      \foreach \direction in {-1,1}
        \draw[trajectory,oriented] plot[domain=0:106,samples=32,variable=\t]
          ({\direction*\t}:{1.05+\side*.44*(tan(\t/2)/tan(53))^2});
    \foreach \angle in {0,180}
      \foreach \tip in {.61,1.49}
         \draw[trajectory,oriented] (\angle:1.05) -- (\angle:\tip);
    \draw[core,oriented] (0:1.05) arc (0:180:1.05);
    \draw[core,oriented] (0:1.05) arc (0:-180:1.05);
    \fill (0:1.05) circle (2pt);
    \fill (180:1.05) circle (2pt);
    \node[inner sep=0pt] at (.80,-.23) {$e$};
    \node[inner sep=0pt] at (-.80,-.23) {$h$};
    \node at (0,-1.82) {Legendrian loop edge};
  \end{scope}
\end{tikzpicture}
\caption{An annular example in Hayden's ribbon construction. Inserting a positive elliptic point $e$ and a positive hyperbolic point $h$ breaks the repelling closed leaf. Both stable branches of $h$ originate at $e$. The two heavy arcs on the right, together with $e$ and $h$, form one Legendrian circle in the ribbon graph.}
\label{fig:ribbon-spine}
\end{figure}
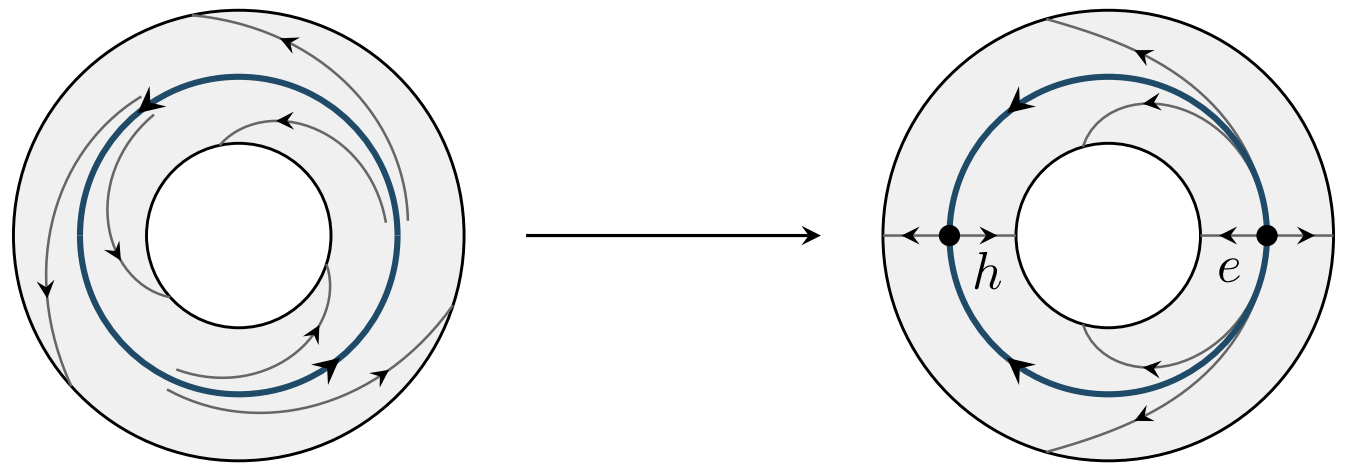

\begin{lemma}[The terminal positive region]\label{lem:terminal}
Suppose a protected configuration has no interior negative annuli, so that
$R_-=\bigsqcup_iN_i$ and $\Gamma_F=\bigsqcup_iC_i$. 
Then $F$ is isotopic relative to its transverse boundary to a Legendrian ribbon. 
\end{lemma}
\begin{proof}
Use the contact vector field $X$ defining the dividing set to give a flow neighborhood $(-\varepsilon,\varepsilon)_t\times F$ of the whole surface.
Since $X$ is a contact vector field, we have
$\Lie_X\alpha=g\alpha$ for some smooth function $g$.
We first rescale $\alpha$ by a positive smooth function
$u$ satisfying
\[
X(\log u)=-g,\qquad u|_F=1.
\]
Since $X$ is transverse to $F$, in the flow coordinates
$(p,t)\in F\times(-\varepsilon,\varepsilon)$ we have
$X=\partial_t$. Thus the above equation reduces to
an ordinary differential equation in $t$ for each
$p\in F$, and its solution is explicitly given by
\[
u(p,t)=\exp\left(-\int_0^t g(p,s)\,ds\right).
\]
Replacing $\alpha$ by $u\alpha$ makes it flow-invariant.
Indeed,
\[
\begin{aligned}
\Lie_X(u\alpha)
&=X(u)\alpha+u\Lie_X\alpha\\
&=u\bigl(X(\log u)+g\bigr)\alpha\\
&=0.
\end{aligned}
\]
Hence, denoting the rescaled contact form again by
$\alpha$, we may write
\[
\alpha=f\,dt+\beta,
\]
where $f$ and $\beta$ are independent of $t$,
$f^{-1}(0)=\Gamma_F$, and $f>0$ on the interiors
of the positive regions.
The contact condition becomes
\[
f\,d\beta+\beta\wedge df>0.
\]

On the interior $P$ of each positive component set
\[
\lambda=\dfrac{\beta}{f},\qquad
\omega=d\lambda
=\frac{f\,d\beta+\beta\wedge df}{f^2}>0. 
\]
Since $\omega$ is nondegenerate, one can define a smooth vector field $Z$ by
\[
\iota_Z\omega=\lambda. 
\]
Then $Z$ generates the oriented characteristic foliation and
\begin{equation}\label{eq:liouvilleflow}
\Lie_Z\omega=d(\iota_Z\omega)=\omega. 
\end{equation}
At a singularity, contactness and $f>0$ give $d\beta>0$, so the singularity is positive. 
By the return-map calculation above, every periodic orbit in $P$ is repelling, with multiplier $e^{T_\gamma}$. 

We truncate $P$ along small positive levels $f=\delta_i$ near its boundary circles $C_i$, obtaining a compact surface $Q$. 
To check the boundary direction, fix a positive area form $\Omega$ on $F$ and define $W$ by $\iota_W\Omega=\beta$. 
On $\Gamma_F$ the contact inequality becomes
\[
 0<\beta\wedge df=-W(f)\Omega. 
\]
Thus $W(f)<0$ on sufficiently small collars of the $C_i$. 
Since $Z$ is a positive multiple of $W$ on $P$, it points strictly out of $Q$. 
The characteristic trajectories in the removed strips continue through the protected product annuli $N_i$ and exit at $T_i$; these strips and annuli contain no singularities or recurrent trajectories. 

Apply Lemma~\ref{lem:positive-genericity} on each $Q$, retaining $Z$ near its boundary, and set $\lambda'=\iota_{Z'}\omega$. 
Then $\lambda'=\lambda$ near $\partial Q$ and $d\lambda'>0$. 
Define a target foliation on the whole abstract surface $F$ by the oriented kernels of
\[
 \beta'=\begin{cases}
 f\lambda'&\text{on each truncated positive region }Q,\\
 \beta&\text{elsewhere}.
 \end{cases}
\]
These forms agree on full collars of $\partial Q$. 
Moreover, 
\[
 f\,d\beta'+\beta'\wedge df=f^2d\lambda'>0
\]
on each $Q$, and the original inequality is unchanged elsewhere. 
In particular, every singularity of the target foliation in $Q$ is positive, since at such a point $p$, 
\[
\lambda'(p)=0,
\qquad
d\beta'(p)=f(p)d\lambda'(p)>0.
\]
Outside $Q$, the target foliation agrees with the original one, whose singularities are positive. 
Consequently the target foliation is divided by the original $\Gamma_F$ and agrees with the original foliation on every component of $F\setminus\Gamma_F$ meeting $\partial F$, namely the protected negative collars. 
The relative realization theorem \cite[Theorem~2.2]{EVHM} therefore gives a convex isotopy of the whole surface, fixed near its boundary, realizing this target foliation. 

The realized foliation is Morse--Smale, has only positive singularities, and has no attracting closed leaves. 
Indeed, all nonwandering points lie in the interiors of
the truncated positive regions $Q$, where the genericity lemma applies. 
Moreover, every periodic orbit is repelling. 
The protected strips and annuli contain no nonwandering points, and trajectories entering them exit through the transverse boundary. Thus the Morse--Smale property of the dynamics in the truncated positive regions gives the Morse--Smale property of the realized foliation on the whole surface. 
Hayden's ribbon-flexibility lemma~\cite[Lemma~3.1]{Hayden} therefore yields an isotopy, relative to the transverse boundary, to a Legendrian ribbon.
In this process, each repelling closed leaf is broken by inserting a positive elliptic--hyperbolic pair, as illustrated in Figure~\ref{fig:ribbon-spine}.
If $F$ is disconnected, the isotopies can be performed in disjoint neighborhoods of its components, yielding a Legendrian ribbon for the disjoint union of the resulting Legendrian graphs. 
\end{proof}

\section{Annulus elimination and the main theorem}
We first state the two local inputs.
Their proofs are given in Sections 4 and 5; the rest of this section explains how they give a finite elimination procedure.
\subsection{The two local inputs}
\begin{proposition}[Relative annulus replacement]\label{prop:L}
Let $(Y, \xi)$ be a closed connected contact three-manifold whose Ozsv\'ath--Szab\'o contact invariant 
$c(\xi)\in\HFhat(-Y;\F)$ is nonzero.
Suppose it contains an embedded solid torus of the following form, which we call the \emph{rotation model}:
\begin{equation}\label{eq:model}
 (V,\xi|_V)=\left(D\times S^1_y,
 \ker(\sin x\,dy+\cos x\,dz)\right),
\end{equation}
where $D$ is a rounded rectangle $[-\varepsilon,2\pi+\varepsilon]\times[-\delta,\delta]$ in the $(x,z)$-plane,
$0<\varepsilon<\pi/4$, and $\delta>0$. Let $A_0=\{z=0\}\cap V$, oriented
by $dx\wedge dy$. The contact vector field $\partial_z$ is positively
transverse to $A_0$.
Then there  are a contact structure $\xi'$ on $Y$, a properly embedded
annulus $A'\subset V$, and a contact vector field $X$ for $\xi'$
defined near $A'$, with the following properties.
\begin{enumerate}
\item\label{replacement:relative}
The contact structure $\xi'$ equals $\xi$ outside $V$ and on a
fixed open boundary collar of $V$.
\item\label{replacement:nonzero}
The closed contact invariant $c(\xi')$ is nonzero.
\item\label{replacement:annulus}
The annulus $A'$ agrees with $A_0$ on a collar of each of its
two boundary circles and is smoothly isotopic to $A_0$
relative to these collars.
\item\label{replacement:positive}
The field $X$ is positively transverse to $A'$.
Along $A'$, it is also positively transverse to $\xi'$.
Consequently its dividing set on $A'$ is empty.
\item\label{replacement:ends}
The field $X$ equals $\partial_z$ on a full three-dimensional
neighborhood in $V$ of each boundary circle of $A'$.
\end{enumerate}
\end{proposition}
Proof of this proposition is in Section 4.2.
\begin{remark}[Why the end vector field is fixed]\label{rem:relative-gluing}
Agreement of the contact structures on an open collar of $\partial V$
allows the new filling to be glued to the unchanged exterior. The additional
requirement $X=\partial_z$ on three-dimensional neighborhoods of the annular
ends is used to glue the replacement annulus to the rest of the convex
surface together with its transverse contact vector field. Indeed,
Proposition~\ref{prop:P} makes the exterior field equal to $\partial_z$
there. Thus the replacement deletes the dividing curves on $A_0$ while
leaving the dividing set elsewhere unchanged.
\end{remark}

\begin{figure}[tbp]
\centering
\begingroup
\definecolor{ocblue}{RGB}{31,75,105}
\definecolor{ocgreen}{RGB}{35,104,78}
\begin{tikzpicture}[
  x=1cm,y=1cm,font=\small,
  actual/.style={draw=black,line width=.9pt},
  selected/.style={draw=black,line width=.85pt,
    dash pattern=on 3pt off 2pt},
  outlet/.style={draw=ocgreen,line width=1pt,
    dash pattern=on 5pt off 2pt},
  cut/.style={draw=ocblue,line width=1pt,
    dash pattern=on .7pt off 1.8pt},
  guide/.style={draw=black!60,line width=.4pt},
  every node/.style={inner sep=1.5pt}
]
\node[anchor=west] at (0,2.85)
  {The outlet condition holds};
\fill[black!15] (0,0) rectangle (.70,1.40);
\fill[ocblue!10] (3.10,0) rectangle (6.30,1.40);
\fill[black!20] (3.65,0) rectangle (5.75,1.40);
\fill[black!15] (8.40,0) rectangle (9.45,1.40);
\fill[black!15] (12.50,0) rectangle (13.20,1.40);
\draw[black!45,line width=.4pt] (0,0)--(13.20,0);
\draw[black!45,line width=.4pt] (0,1.40)--(13.20,1.40);
\foreach \x in {0,13.20}
  \draw[actual] (\x,0)--(\x,1.40);
\foreach \x in {.70,8.40,9.45,12.50}
  \draw[outlet] (\x,0)--(\x,1.40);
\foreach \x in {3.65,5.75}
  \draw[selected] (\x,0)--(\x,1.40);
\foreach \x in {3.10,6.30}
  \draw[cut] (\x,0)--(\x,1.40);
\foreach \y in {0,1.40}
  \draw[->,>=stealth,line width=.65pt]
    (10.35,\y)--(11.45,\y);
\node[font=\footnotesize,align=center] at (10.90,2.04)
  {top and bottom edges\\identified};
\node at (.35,.70) {$N_0$};
\node at (1.90,.70) {$P_0$};
\node at (4.70,.82) {$A$};
\node[font=\scriptsize] at (4.70,.40) {selected};
\node at (7.35,.70) {$P_1$};
\node at (8.925,.70) {$A_j$};
\node at (10.90,.70) {$P_2$};
\node at (12.85,.70) {$N_1$};
\node at (0,1.76) {$T_0$};
\node at (13.20,1.76) {$T_1$};
\node[text=ocgreen] at (.70,-.39) {$C_0$};
\node[text=ocgreen] at (8.40,-.39) {$b$};
\node[text=ocgreen] at (9.45,-.39) {$b'$};
\node[text=ocgreen] at (12.50,-.39) {$C_1$};
\draw[ocblue,line width=.55pt]
  (3.10,1.77)--(3.10,2.02)--(6.30,2.02)--(6.30,1.77);
\node[text=ocblue] at (4.70,2.29) {$U$};
\node[text=ocblue,font=\footnotesize] at (3.10,-.39) {$\partial U$};
\node[text=ocblue,font=\footnotesize] at (6.30,-.39) {$\partial U$};
\draw[guide] (3.65,-.08)--(3.65,-.73)--(5.75,-.73)--(5.75,-.08);
\node[fill=white] at (4.70,-.73) {$\partial A$};
\node at (6.60,-1.30)
  {Outlets: $P_0:\ C_0$,\qquad $P_1:\ b$,\qquad
    $P_2:\ b',C_1$.};

\draw[outlet] (.15,-2.12)--(.86,-2.12);
\node[anchor=west] at (1.02,-2.12)
  {unselected dividing curves (outlets)};
\draw[selected] (7.80,-2.12)--(8.51,-2.12);
\node[anchor=west] at (8.67,-2.12) {$\partial A$};
\draw[cut] (10.48,-2.12)--(11.19,-2.12);
\node[anchor=west] at (11.35,-2.12) {$\partial U$};
\fill[black!17] (.15,-2.76) rectangle (.86,-2.52);
\node[anchor=west] at (1.02,-2.64) {negative annuli};
\fill[ocblue!10] (5.18,-2.76) rectangle (5.89,-2.52);
\node[anchor=west] at (6.05,-2.64) {positive collars in $U$};
\end{tikzpicture}
\endgroup
\caption{Outlets for a selected interior negative annulus $A$.
The annulus $U$ contains $A$ and small collars in the adjacent
positive regions. Its dotted boundary curves are not dividing curves
and are not outlets. The top and bottom edges are identified,
so each vertical line represents a circle. Every remaining positive
component has an outlet: $C_0,C_1\subset C$ and
$b,b'\subset\partial A_j$, where $A_j$ is unselected.}
\label{fig:outlet-condition}
\end{figure}
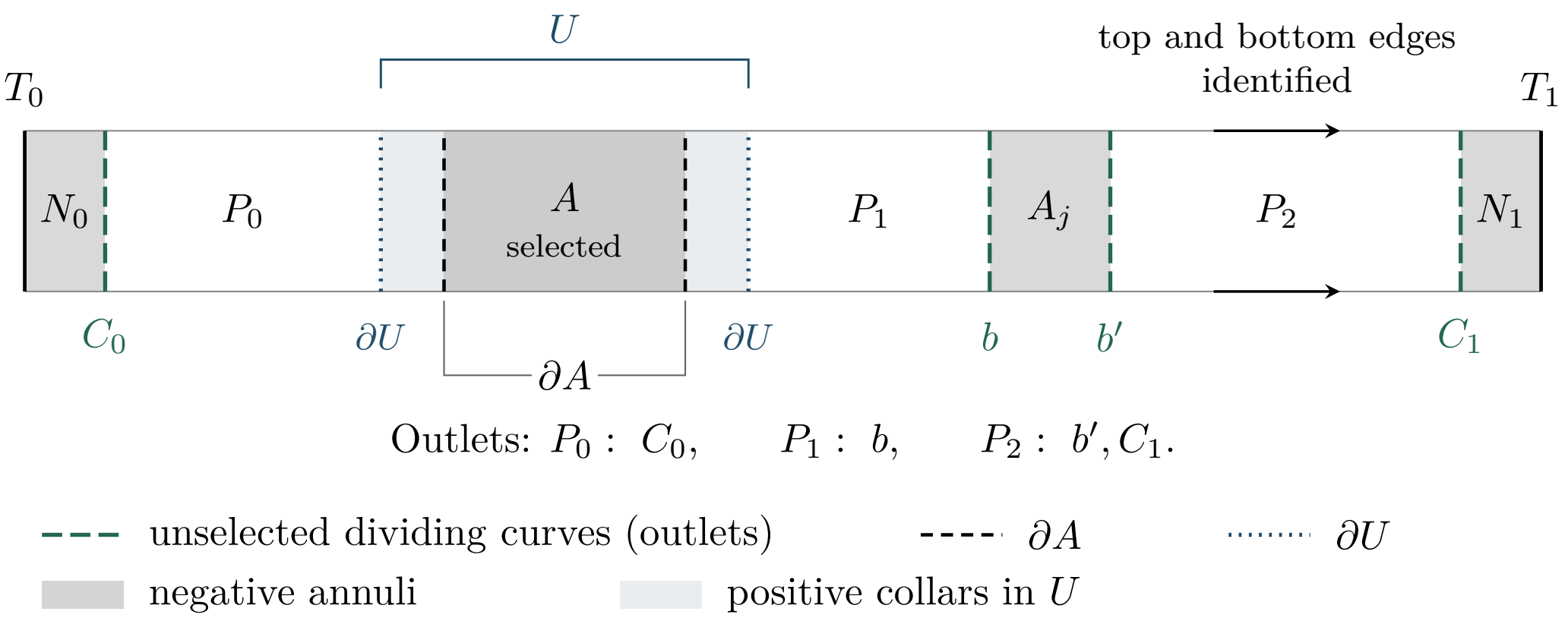

\begin{definition}[Outlets and the outlet condition]\label{def:outlet}
Let $A$ be a selected interior negative annulus in a protected
configuration. Choose an annulus $U$ containing $A$ and small collars
of its two boundary curves in the adjacent positive regions. A component
of $\Gamma_F\setminus\partial A$ is called an \emph{unselected dividing
curve}. If $P$ is a component of
$\overline{R_+(F)\setminus U}$, an \emph{outlet of $P$} is a boundary
component of $P$ that is an unselected dividing curve. Thus an outlet
lies either in the protected multicurve $C$ or in the boundary of an
unselected interior negative annulus. We say that $U$ satisfies the
\emph{outlet condition} if every such $P$ has an outlet. See Figure~\ref{fig:outlet-condition}.
\end{definition}


\begin{proposition}[Relative preparation]\label{prop:P}
Let $(Y,\xi)$ be a positive cooriented closed contact three-manifold.
Let $i_0:F\hookrightarrow(Y,\xi)$ be a protected annular configuration,
and let $X_0$ be the chosen positively transverse contact vector field.
Write $\Gamma_0,R_{+,0},R_{-,0}\subset F$ for its dividing set and
closed signed regions, pulled back by $i_0$.
Let $A$ be an interior negative annulus whose neighborhood
$U\subset F$ satisfies the outlet condition.
Then there are an isotopy of embeddings $i_s$, $0\leq s\leq1$,
and a positively transverse contact vector field $X_1$ near $i_1(F)$
with the following properties.
\begin{enumerate}
\item\label{prep:boundary}
The isotopy is fixed on a neighborhood of $\partial F$.
The characteristic foliation of $i_1(F)$, pulled back by $i_1$,
agrees on the protected collars $N$ with the original foliation.
\item\label{prep:regions}
Let $\Gamma_1,R_{+,1},R_{-,1}\subset i_1(F)$ be the dividing set
and closed signed regions defined using $X_1$. Then
\[
 i_1^{-1}(\Gamma_1)=\Gamma_0,\qquad
 i_1^{-1}(R_{+,1})=R_{+,0},\qquad
 i_1^{-1}(R_{-,1})=R_{-,0}.
\]
These are equalities of subsets of the parameter surface $F$.
\item\label{prep:model}
A neighborhood of $i_1(U)$ has the rotation model
\eqref{eq:model}, with $i_1(U)$ identified with the annular
patch $\{z=0\}$.
\item\label{prep:field}
The field $X_1$ equals $\partial_z$ on this model neighborhood
and equals $X_0$ on a neighborhood of $T=i_0(\partial F)$.
\end{enumerate}
\end{proposition}
Proof of this proposition is in Section 5.3.

\subsection{Finite elimination}

We now recall the folding operation for convex surfaces,
which will be used in the proof of
Theorem~\ref{thm:elimination}.

Let $(Y,\xi)$ be a positive cooriented contact
three-manifold, and let $F\subset Y$ be a compact
oriented convex surface, possibly with boundary.
Choose a contact vector field $X$ positively
transverse to $F$, and denote its dividing set
and positive and negative regions by
$\Gamma_F$, $R_+(F)$, and $R_-(F)$, respectively.

Let $\gamma\subset\operatorname{Int}R_+(F)$ be
a nonisolating simple closed curve. By the
Legendrian realization principle
\cite[Theorem~3.7]{Honda}, we may isotope $F$
so that $\gamma$ becomes Legendrian while preserving
the dividing set. Since $\gamma$ is disjoint
from $\Gamma_F$, its contact twisting relative
to $F$ is
\[
  \operatorname{tw}(\gamma,F)
  =-\frac12\#(\gamma\cap\Gamma_F)=0.
\]

We can now perform a \emph{fold along $\gamma$}.
This is a local isotopy of $F$, supported in a
small neighborhood of $\gamma$, which produces
a new convex surface with two additional dividing
curves parallel to $\gamma$. We choose the fold
so that these two curves bound a new negative
annulus $A_{\mathrm{aux}}$. Outside the annular
neighborhood where the fold is performed, the
dividing set and the signed regions remain
unchanged, under the identification induced
by the isotopy.

\begin{figure}[tbp]
\centering
\begin{tikzpicture}[x=1cm,y=1cm,font=\small,
  seam/.style={draw=black!45,densely dotted,line width=.5pt},
  edge/.style={draw=black!45,line width=.4pt},
  divider/.style={draw=black,line width=.85pt}]
  \foreach \y in {-1,1}
    \draw[edge] (0,\y) -- (3.7,\y);
  \foreach \x in {0,3.7}
    \draw[seam] (\x,-1) -- (\x,1);
  \draw[dashed,figureblue,line width=.7pt] (0,0) -- (3.7,0)
    node[right=2pt,text=black] {$\gamma$};
  \node at (1.85,.6) {$R_+$};
  \node at (1.85,-.6) {$R_+$};
  \draw[->,black!70] (4.4,0) -- node[above=3pt] {fold} (5.7,0);
  \begin{scope}[xshift=6.25cm]
    \fill[black!9] (0,-.38) rectangle (3.7,.38);
    \foreach \y in {-1,1}
      \draw[edge] (0,\y) -- (3.7,\y);
    \foreach \x in {0,3.7}
      \draw[seam] (\x,-1) -- (\x,1);
    \foreach \y in {-.38,.38}
      \draw[divider] (0,\y) -- (3.7,\y);
    \node at (1.85,.72) {$R_+$};
    \node at (1.85,-.72) {$R_+$};
    \node at (1.85,0) {$A_{\mathrm{aux}}\subset R_-$};
  \end{scope}
\end{tikzpicture}
\caption{The change of dividing set under a fold in a positive region.
The vertical sides of each strip are identified. Before folding, the
dashed core $\gamma$ is not a dividing curve. Afterwards, the two
bold curves bound the new negative annulus $A_{\mathrm{aux}}$.
The thin horizontal edges delimit the chosen annular neighborhood.
The diagram shows the signed regions, rather than the spatial shape
of the folding isotopy.}
\label{fig:fold-local}
\end{figure}
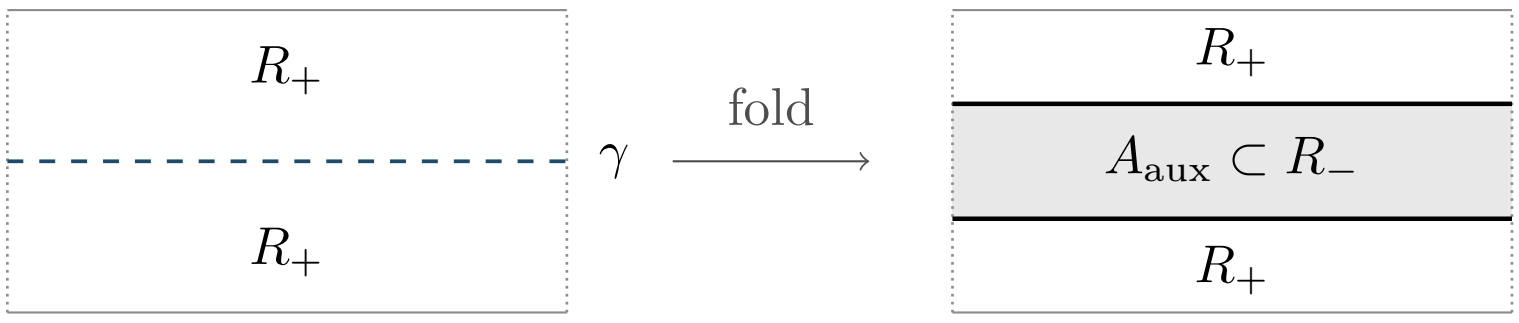

The ambient contact structure $\xi$ is not
modified by this operation. Although the initial
and final surfaces are convex, the intermediate
surfaces need not be convex. This is Honda's
folding construction
\cite[\S\S3.3.1 and 5.3.1, Figure~19]{Honda}.
Figure~\ref{fig:fold-local} illustrates the
change in the dividing set and the signed regions.

\begin{theorem}[Annulus elimination]\label{thm:elimination}
Let $(F,\xi)$ be a protected annular configuration in a closed connected contact
three-manifold with $c(\xi)\ne0$. For any chosen interior negative annulus $A$,
a finite operation produces a new contact structure and a smoothly
boundary-relative isotopic surface such that the closed contact invariant
remains nonzero and
\begin{equation}\label{eq:delete}
 \Gamma_{\mathrm{new}}=\Gamma_{\mathrm{old}}\setminus\partial A.
\end{equation}
The dividing sets are compared using the surface isotopy. The operation
uses at most one fold and two local replacements. The boundary is fixed
near $T$, and the pulled-back protected-collar foliation is unchanged.
If initially there are $m_0$ interior negative annuli, at most $2m_0$
local replacements eliminate all of them.
\end{theorem}

\begin{proof}
Let $m$ denote the current total number of interior negative annuli,
counted over all components of $F$. Define a graph $G$ with one vertex
for each connected component of $R_+(F)$ and one edge for each interior
negative annulus. An annulus with boundary curves $a_0,a_1$ gives an
edge whose endpoints are the positive regions incident to $a_0$ and
$a_1$, respectively. If these regions coincide, the edge is a loop;
its two ends record two distinct boundary components of that positive
region. Multiple edges are allowed, and a loop contributes two to the
degree of its vertex. Mark a vertex when its positive region has a
boundary component $C_i$.

The surface $F$ is obtained by joining its positive regions with the
interior negative annuli and then attaching the protected collars.
Each protected collar attaches to just one positive region and hence
does not join two components. Consequently two positive regions lie
in the same component of $F$ exactly when their vertices are joined
by a path in $G$. Every component of $F$ has nonempty boundary, so
every component of $G$ contains a marked vertex.

If the selected annulus has the outlet condition, prepare it by
Proposition~\ref{prop:P} and replace it by Proposition~\ref{prop:L}.
The latter proposition supplies a positive transverse contact vector
field agreeing with the exterior field on open end neighborhoods.
The replacement therefore glues to the convex remainder and deletes
exactly $\partial A$. The adjacent positive regions merge if they were
distinct; otherwise their connectivity is unchanged. In either case
the number of interior negative annuli decreases by one.

Suppose the outlet condition fails. Removing $U$ removes only boundary
collars from positive components, so any additional incident edge or
protected boundary component supplies an outlet. A loop with no other
incident edge and no protected boundary would form a graph component
with no marked vertex, which is impossible. Thus a failing side is an
unmarked vertex $P$ of degree one. Both ends of $A$ cannot be such
vertices, since they would again form a component with no marked vertex.
There is exactly one failing side.

The positive region $P$ has one boundary curve and cannot be a disk:
otherwise the core of $A$ would bound a disk in $F$, contrary to
Lemma~\ref{lem:annular-core}. Here $c(\xi)\ne0$ supplies the required
tightness. Hence $P$ has positive genus. Choose a nonseparating simple
closed curve $\gamma\subset\Int P$.
Cutting along $\gamma$ leaves $P$ connected. Moreover, $\gamma$ is
disjoint from $\partial P$, so this connected complement still has
the original boundary component $\partial P\subset\Gamma_F$.
All other components of $F\setminus\Gamma_F$ are unchanged.
Thus every component of $F\setminus(\Gamma_F\cup\gamma)$ has closure
meeting $\Gamma_F$, which is the nonisolating condition. Relative Legendrian
realization makes it Legendrian while retaining the protected boundary
foliation \cite[Theorem~3.7]{Honda}, \cite[Theorem~2.2]{EVHM}.
Apply the fold just described along $\gamma$, creating
$A_{\mathrm{aux}}$. Since $\gamma$ is nonseparating in $P$, the
remaining positive part of $P$ is connected, and its two new dividing
boundaries are outlets for the selected annulus $A$.

Prepare and replace $A$, choosing its replacement torus after preparation
and disjoint from $A_{\mathrm{aux}}$. The resulting dividing set is
\[
 (\Gamma_{\mathrm{old}}\setminus\partial A)
 \sqcup\partial A_{\mathrm{aux}}.
\]
The former positive leaf has merged with its neighbor, so
$A_{\mathrm{aux}}$ is now a loop at a single positive component $Q$.
If $Q$ is marked, some $C_i$ supplies an outlet. Otherwise a path from
$Q$ to a marked vertex starts with an edge distinct from this loop and
supplies an outlet. Thus $A_{\mathrm{aux}}$ itself satisfies the outlet
condition. Prepare and replace it to remove the two auxiliary curves,
proving \eqref{eq:delete}. Figure~\ref{fig:fold-elimination} shows these
graph changes when the neighbor of $P$ is marked.

\begin{figure}[tbp]
\centering
\begin{tikzpicture}[
  x=1cm, y=1cm, font=\small,
  region/.style={circle,draw=black,fill=white,line width=.8pt,
                 minimum size=5.2pt,inner sep=0pt},
  marked/.style={region,fill=black},
  annulus/.style={draw=black,line width=.85pt},
  auxiliary/.style={draw=figureblue,line width=1pt},
  operation/.style={->,draw=black!72,line width=.55pt},
  every node/.append style={outer sep=0pt}
]
\draw[annulus] (0,0) -- node[above=2pt] {$A$} (1.35,0);
\node[marked] at (0,0) {};
\node[region] at (1.35,0) {};
\node[below=3pt] at (1.35,0) {$P$};
\draw[annulus] (3.5,0) -- node[above=2pt] {$A$} (4.85,0);
\draw[auxiliary] (4.85,0)
  .. controls (4.40,.33) and (4.40,.82) .. (4.85,.82)
  .. controls (5.30,.82) and (5.30,.33) .. (4.85,0);
\node[above=2pt] at (4.85,.82) {$A_{\mathrm{aux}}$};
\node[marked] at (3.5,0) {};
\node[region] at (4.85,0) {};
\draw[auxiliary] (8.20,0)
  .. controls (7.75,.33) and (7.75,.82) .. (8.20,.82)
  .. controls (8.65,.82) and (8.65,.33) .. (8.20,0);
\node[marked] at (8.20,0) {};
\node[below=3pt] at (8.20,0) {$Q$};
\node[marked] at (11.9,0) {};
\draw[operation] (1.72,0) -- node[above=2pt] {fold} (3.12,0);
\draw[operation] (5.49,0) -- node[above=2pt] {replace $A$} (7.53,0);
\draw[operation] (8.87,0) -- node[above=2pt] {replace $A_{\mathrm{aux}}$} (11.45,0);
\foreach \x/\counttext in {.675/m,4.175/{m+1},8.20/m,11.9/{m-1}}
  \node at (\x,-.83) {$\counttext$};
\end{tikzpicture}
\caption{The case in which the outlet condition fails, shown when the
unmarked leaf $P$ has a marked neighbor (solid vertex). Folding along a nonseparating
curve in $P$ creates the loop $A_{\mathrm{aux}}$. Replacing $A$ merges the
two positive regions. A protected boundary of the merged region $Q$ then
supplies an outlet for replacing $A_{\mathrm{aux}}$. The bottom row counts
all interior negative annuli.}
\label{fig:fold-elimination}
\end{figure}
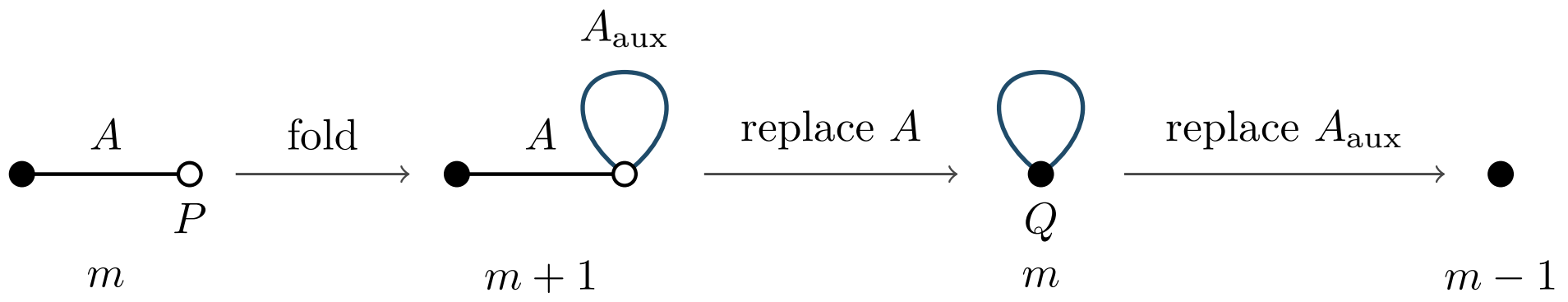

Each application of Proposition~\ref{prop:L} preserves nonvanishing
of the closed contact invariant. Folding and preparation are surface
isotopies and do not change the ambient contact structure. When the
outlet condition initially fails, the total number $m$ of interior
negative annuli changes as follows:
\[
 m\ \xrightarrow{\mathrm{fold}}\ m+1
 \ \xrightarrow{\mathrm{replace}\ A}\ m
 \ \xrightarrow{\mathrm{replace}\ A_{\mathrm{aux}}}\ m-1.
\]
Each completed operation decreases this total by one, uses at most
two applications of Proposition~\ref{prop:L}, and preserves the
property that every graph component contains a marked vertex.
This proves the bound $2m_0$.
Every replacement fixes the surface outside its selected solid torus,
so the argument applies to all connected components simultaneously.
\end{proof}

\subsection{Proof of the main theorems}
\begin{proof}[Proof of Theorem~\ref{main:surface}]
Apply Proposition~\ref{prop:normal-form} to the whole surface $F$.
Initially $c(\xi_{\std})\ne0$ \cite{OSContact}. Theorem~\ref{thm:elimination} eliminates
the interior negative annuli in finitely many operations, with the total
count taken over every connected component. All changes to $F$ are
smooth isotopies fixed near its boundary. Each local contact replacement
preserves a nonzero closed contact invariant and fixes the surface outside
its replacement torus. Thus the entire resulting surface lies in one
contact structure $\xi_{\mathrm{fin}}$ with
$c(\xi_{\mathrm{fin}})\ne0$.

By Lemma~\ref{lem:terminal}, this surface is isotopic to a Legendrian
ribbon in $(S^3,\xi_{\mathrm{fin}})$. Nonvanishing of the contact invariant
implies that $\xi_{\mathrm{fin}}$ is tight \cite{OSContact}. The positive tight contact
structure on $S^3$ is unique up to isotopy \cite{Eliashberg}. Gray stability
therefore gives a diffeomorphism isotopic to the identity carrying
$\xi_{\mathrm{fin}}$ to $\xi_{\std}$. Apply this one diffeomorphism to
the whole ribbon. The ribbon-to-quasipositive-surface correspondence
\cite{BI,Hayden} identifies its smooth ambient isotopy class with that
of a quasipositive surface. This proves the assertion for the entire
embedded surface, including when it is disconnected.
\end{proof}

\begin{proof}[Proof of Theorem~\ref{main:link}]
Suppose $\slbar(L)=-\chi_{\max}(L)$. Choose a transverse representative
$T$ with $\slink(T)=\slbar(L)$, and transport any Seifert surface $F$
realizing $\chi_{\max}(L)$ to one with boundary $T$ by a smooth ambient
isotopy. Since self-linking in $S^3$ is independent of the spanning
surface,
\[
 \slink_F(T)=\slbar(L)=-\chi_{\max}(L)=-\chi(F).
\]
Theorem~\ref{main:surface} makes $F$ quasipositive up to ambient isotopy.
This proves both strong quasipositivity of $L$ and the assertion about
every surface realizing $\chi_{\max}(L)$.

Conversely, let $F_\beta$ be a positive band surface for a strongly
quasipositive braid representative of $L$. Its transverse closure has
self-linking $-\chi(F_\beta)$. Therefore
\[
 -\chi(F_\beta)\le\slbar(L)\le-\chi_{\max}(L)
 \le-\chi(F_\beta).
\]
All inequalities are equalities, which proves the converse.
\end{proof}

\begin{remark}[The transverse boundary]\label{rem:transverse-boundary}
The preparation and elimination isotopies fix the boundary $T$, and
every contact replacement is supported away from $T$. Since there are
finitely many replacements, the final contact structure
$\xi_{\mathrm{fin}}$ agrees with $\xi_{\std}$ on some fixed open
neighborhood $N(T)$. Lemma~\ref{lem:terminal} then produces a
Legendrian ribbon $R$ in $(S^3,\xi_{\mathrm{fin}})$ with
$\partial R=T$, smoothly ambiently isotopic to $F$ relative to $T$.

The remaining issue is the identification of $\xi_{\mathrm{fin}}$
with $\xi_{\std}$. Tightness and the uniqueness theorem for the
standard contact sphere provide an isotopy of contact structures on
$S^3$ \cite[Theorem~2.1.1]{Eliashberg}. They do not supply a contact
path fixed on $N(T)$. For a path $\xi_s$ from $\xi_{\mathrm{fin}}$
to $\xi_{\std}$, Gray stability gives an ambient isotopy $g_s$,
with $g_0=\mathrm{id}$ and $g_s^*\xi_s=\xi_{\mathrm{fin}}$.
Thus $g_s(T)$ is positively transverse to $\xi_s$. Its transversality
to the fixed contact structure $\xi_{\std}$ at intermediate times
is not established. Consequently, the argument does not identify the
transverse isotopy classes of $T$ and $g_1(T)$ in
$(S^3,\xi_{\std})$. Equality of the endpoint contact structures near
$T$ does not supply the missing relative path in this proof.

If the contact path could be chosen constant on $N(T)$, relative
Gray stability would allow $g_s$ to be the identity there. Then
$g_1(R)$ would be a Legendrian ribbon in the original contact
structure with boundary exactly $T$, and the resulting surface
isotopy would fix $T$. Such a relative path is not constructed here.
The theorem therefore concerns the smooth ambient isotopy class of
the surface. Whether the initially specified transverse isotopy
class of $T$ can be preserved remains unresolved by this argument.
\end{remark}
\section{The six-suture replacement}
We prove Proposition~\ref{prop:L}. Four elementary lemmas will control
transverse contact vector fields, contact forms on boundary collars,
and the replacement annulus.

\subsection{Contact Hamiltonians and boundary collars}
The following Hamiltonian identity will be used in this section and
in Section~5.
\begin{lemma}[Hamiltonian transversality]\label{lem:ham}
Let $(M,\ker\alpha)$ be a positive cooriented contact three-manifold.
Let $i:F\hookrightarrow(M,\ker\alpha)$ be an oriented embedded surface,
and put $\beta=i^*\alpha$. For a smooth function $H$ near $F$, the contact
vector field $X_H$ is characterized by
\[
 \alpha(X_H)=H,\qquad
 \iota_{X_H}d\alpha=(R_\alpha H)\alpha-dH,
\]
where $R_\alpha$ is the Reeb vector field. It satisfies
\begin{equation}\label{eq:ham}
 i^*\bigl(\iota_{X_H}(\alpha\wedge d\alpha)\bigr)
 =(H\circ i)\,d\beta+\beta\wedge d(H\circ i).
\end{equation}
If this two-form is positive, then $X_H$ is positively transverse to $F$
and determines the dividing set $\{H\circ i=0\}$ on the parameter
surface. In particular, the transversality test depends only on
$H\circ i$ and its differential on $F$.
\end{lemma}
\begin{proof}
Writing $X_H=HR_\alpha+V$ with $V\in\ker\alpha$ determines $V$ uniquely
by the nondegeneracy of $d\alpha$ on $\ker\alpha$. Cartan's formula gives
$\Lie_{X_H}\alpha=(R_\alpha H)\alpha$, so its flow preserves the
cooriented contact structure. Contracting the contact volume form gives
\[
 \iota_{X_H}(\alpha\wedge d\alpha)
 =H\,d\alpha-\alpha\wedge\iota_{X_H}d\alpha
 =H\,d\alpha+\alpha\wedge dH.
\]
Its pullback is positive exactly when $X_H$ points in the positive normal
direction. The dividing set is where $\alpha(X_H)=0$.
If two Hamiltonians have the same restriction to $F$, their restrictions
also have the same differential on $F$, so the displayed test is
unchanged. If they agree on a full ambient neighborhood, the defining
equations give the same contact vector field on that neighborhood.
\end{proof}

\begin{lemma}[A contact collar with fixed end forms]\label{lem:collar}
Consider the torus collar $[0,1]_r\times S^1_u\times S^1_y$, oriented by
$dr\wedge du\wedge dy$. Here $r$ is the collar parameter and $u,y$ are
angular coordinates. Let $H,f,b$ be smooth functions of $u$ such that
\[
 W=Hf'-fH'>0,
\]
where the primes denote differentiation with respect to $u$.
Let $k:[0,1]\to[0,1]$ be smooth, monotone, constant near both
endpoints, and satisfy either $(k(0),k(1))=(1,0)$ or
$(k(0),k(1))=(0,1)$. Choose a constant $M$ and define $h$ by
\[
 M>\max_{u\in S^1}\frac{|f(u)b(u)|}{W(u)},\qquad
 h(r)=1+M|k'(r)|.
\]
Then the one-form
\begin{equation}\label{eq:collar}
 \alpha_k=h(r)H(u)\,dr+f(u)\,dy+k(r)b(u)\,du
\end{equation}
is a positive contact form on the collar. For $(k(0),k(1))=(1,0)$,
it equals $H\,dr+f\,dy+b\,du$ on a neighborhood of $r=0$
and $H\,dr+f\,dy$ on a neighborhood of $r=1$.
For $(k(0),k(1))=(0,1)$, these two expressions are reversed. The first choice
will be used for the normalization collar $L_{\rm in}$, and the second
for the restoring collar $L_{\rm out}$ in Proposition~\ref{prop:L}.
\end{lemma}
\begin{proof}
Since $H,f,b$ depend only on $u$ and $h,k$ depend only on $r$,
\[
 d\alpha_k=(k'b-hH')\,dr\wedge du+f'\,du\wedge dy.
\]
Wedging with $\alpha_k$ in the orientation $dr\wedge du\wedge dy$
gives
\begin{align*}
 \alpha_k\wedge d\alpha_k
 &=\bigl(hHf'+f(k'b-hH')\bigr)\,dr\wedge du\wedge dy\\
 &=\bigl(h(Hf'-fH')+fk'b\bigr)\,dr\wedge du\wedge dy\\
 &=\bigl((1+M|k'|)W+fk'b\bigr)\,dr\wedge du\wedge dy.
\end{align*}
Now $fk'b\geq-|k'||fb|$, so its coefficient is at least
\[
 W+|k'|(MW-|fb|)>0,
\]
where the last inequality follows pointwise from the strict choice of $M$
and $W>0$. Monotonicity also makes $|k'|$ smooth: it equals $k'$ for
the increasing choice and $-k'$ for the decreasing choice. On the endpoint
collars $k'=0$, hence $h=1$; inserting the endpoint values of $k$
gives the asserted contact forms on the two end collars.
\end{proof}

\begin{lemma}[Interpolation with a fixed signed zero set]\label{lem:fixedf}
Let $D$ be a compact disk. On $D\times S^1_y$, let
\[
 \alpha_j=f_j\,dy+\beta_j,\qquad j=0,1,
\]
be $S^1$-invariant positive contact forms that agree in a boundary collar.
Assume the $f_j$ have the same signed regular zero set: their zero sets
coincide, zero is a regular value of each $f_j$, and their signs agree on
the complement. Then, after division of $\alpha_1$ by a positive function
equal to one near the boundary, the two forms are connected by an
$S^1$-invariant contact path fixed in that collar.
\end{lemma}
\begin{proof}
The quotient $q=f_1/f_0$ extends to a smooth positive function: in a local
defining coordinate for a regular zero arc both numerator and denominator
are that coordinate times a nonvanishing function of the same sign. Set
\[
 \alpha_t=f_0\,dy+(1-t)\beta_0+t\beta_1/q.
\]
For $f\,dy+\beta$, contactness is equivalent to positivity of
$f\,d\beta+\beta\wedge df$ on the appropriately oriented disk. For fixed
$f=f_0$ this is linear in $\beta$. Thus the path is contact and is constant
in the boundary collar. Relative Gray stability on the compact solid torus
provides the corresponding boundary-relative contact identification.
\end{proof}

We now apply Lemma~\ref{lem:ham} to a product annulus.

\begin{lemma}[A positive convex product annulus]\label{lem:annulus}
Let $A=I_s\times S^1_y$ be embedded in a contact manifold and suppose
\[
 \beta=\alpha|_{TA}=p(s)\,dy+q(s)\,ds.
\]
Assume $p$ has exactly one simple zero, at $0$, and changes from negative to
positive there. Let $r$ be smooth with $r(0)=0$ and $r'>0$. Then $H=p/r$
extends smoothly and positively across $0$, and
\begin{equation}\label{eq:annulusH}
 H\,d\beta+\beta\wedge dH=H^2r'\,ds\wedge dy>0.
\end{equation}
Consequently $A$ admits a positive transverse contact vector field with
empty dividing set.
\end{lemma}
\begin{proof}
The quotient is smooth by the simple-zero assumption and has value
$p'(0)/r'(0)>0$ at $0$. Since $H$ depends only on $s$,
\[
 H\,d\beta+\beta\wedge dH=(Hp'-pH')\,ds\wedge dy.
\]
Substitution of $p=Hr$ gives \eqref{eq:annulusH}. The $q$-term contributes
nothing. Lemma~\ref{lem:ham} completes the argument.
\end{proof}

\begin{remark}\label{rem:endHamiltonian}
If the ends lie in the rotation model
$\alpha=\sin x\,dy+\cos x\,dz$, choose $r=\tan x$ there. Then $H=\cos x$.
Extending $H$ to be $\cos x$ on the full three-dimensional end neighborhoods
makes its contact vector field exactly $\partial_z$ on those neighborhoods.
This is stronger than agreement of a Hamiltonian only along a boundary
circle. The extension elsewhere is unconstrained by
Lemma~\ref{lem:annulus}, apart from smoothness.
\end{remark}

\subsection{Proof of the local replacement}

\begin{proof}[Proof of Proposition~\ref{prop:L}]
We construct two contact structures $\xi_B,\xi_C$ satisfying
\textup{(\ref{replacement:relative})}. For each $K=B,C$ we then
construct an annulus $A'_K$ and a contact vector field $X_K$ satisfying
\textup{(\ref{replacement:annulus})},
\textup{(\ref{replacement:positive})}, and
\textup{(\ref{replacement:ends})}. Finally, the contact-element
relation allows us to choose one index $K$ for which
\textup{(\ref{replacement:nonzero})} also holds. All boundary
identifications below are product maps with respect to the $S^1_y$ factor.

\subsection*{Step 1. The marked boundary}
Write $\alpha_E=\sin x\,dy+\cos x\,dz$. Its contact volume form is
$dx\wedge dy\wedge dz$. The vector field
\[
 X_0=\sin x\cos x\,\partial_x+z\,\partial_z
\]
satisfies $\alpha_E(X_0)=z\cos x$ and
$\Lie_{X_0}\alpha_E=\cos^2x\,\alpha_E$. 
This contact vector field points outward on all four sides of the
rectangle and remains transverse after a sufficiently small rounding.
The dividing set on $\partial D\times S^1_y$ is therefore
$\{z\cos x=0\}\cap\partial V$. On $\partial D$, the equation $z=0$
has two solutions, one on each vertical side. The equation $\cos x=0$
has two solutions on each horizontal side. Taking the product of these
six points with $S^1_y$ gives the six longitudinal dividing circles.

On the meridian disk $D\times\{y_0\}$, the contact vector field
$\partial_y$ is transverse to the disk. Its dividing set is given by
$\alpha_E(\partial_y)=\sin x=0$, namely the three arcs at
$x=0,\pi,2\pi$.

Number their endpoints clockwise, with $1,2,3$ on the upper side from
left to right, and $4,5,6$ on the lower side from right to left. The
three arcs join $1$ to $6$, $2$ to $5$, and $3$ to $4$.
We denote this pairing of the six endpoints by
\[
 E=(16)(25)(34).
\]
Let $P=(-\varepsilon,0)$ and $Q=(2\pi+\varepsilon,0)$ on the marked disk
boundary. They lie in the intervals $(6,1)$ and $(3,4)$, respectively.
The $dy$ coefficient of $\alpha_E$ is negative at $P$ and positive at $Q$.
Fix small product neighborhoods of these points throughout the proof.

Let $\phi_r$ be the flow of $X_0$, and use
$\Psi_0(r,p)=\phi_r(p)$ as collar coordinates, with $r=0$ on
$\partial V$ and $r$ increasing outwards. Extend a clockwise coordinate
$u$ on $\partial D$ and the coordinate $y$ constantly along this flow.
Since $\Lie_{X_0}\alpha_E=\cos^2x\,\alpha_E$, the rescaling
\[
 \alpha_{\rm old}=
 \exp\!\left(-\int_0^r\cos^2\!\bigl(x(\phi_s(p))\bigr)\,ds\right)
 \Psi_0^*\alpha_E
\]
is independent of $r$; it is also $S^1_y$-invariant. Hence
\begin{equation}\label{eq:old-collar-form}
 \alpha_{\rm old}=H(u)\,dr+f(u)\,dy+b(u)\,du,
 \qquad W=Hf'-fH'>0.
\end{equation}
The inequality follows from contactness in the orientation
$dr\wedge du\wedge dy$. Evaluating at $r=0$ gives
\[
 H=z\cos x,\qquad f=\sin x,\qquad b=\cos x\,\frac{dz}{du},
 \qquad (x,z)=p(u)\in\partial D.
\]
Thus $f$ has six simple zeros, the endpoints of the meridional
dividing arcs, while the zeros of $H$ give the six torus sutures.
These two sets alternate around $\partial D$.

\subsection*{Step 2. Collar normalization and three fillings}
Let $L_{\rm in}=[0,1]_r\times\partial V$ be a \emph{normalization
collar}. Glue its $r=0$ end to $\partial V$ in the outward flow
coordinates from Step~1. On $L_{\rm in}$ use the form $\alpha_k$
of \eqref{eq:collar}, with $H,f,b$ from \eqref{eq:old-collar-form} and
$k(0)=1$, $k(1)=0$.
Near $r=0$ the function $k$ is identically $1$ and $h=1$, so the
form agrees with $\alpha_{\rm old}$ in the flow coordinates and extends
it smoothly across the gluing torus. Near $r=1$ it equals $H\,dr+f\,dy$.
The coefficient $\alpha_k(\partial_y)=f(u)$ is independent of
$r$, so each meridional dividing arc just acquires a radial endpoint
extension. The matching remains $E$ with its signs.

Since $W>0$, the pair $(H,f)$ never vanishes simultaneously. Write
$H+if=R e^{i\theta}$, where $R=\sqrt{H^2+f^2}>0$ and $\theta$ is a
continuous real lift of the argument along the oriented boundary.
Differentiating this identity with respect to $u$ gives
\[
 H'+if'=(R'+iR\theta')e^{i\theta},\qquad
 \frac{H'+if'}{H+if}=\frac{R'}{R}+i\theta'.
\]
Taking imaginary parts, or multiplying the quotient by
$(H-if)/(H-if)$, yields
\begin{equation}\label{eq:angle-derivative}
 \theta'=\operatorname{Im}\!\left(\frac{H'+if'}{H+if}\right)
 =\frac{Hf'-fH'}{H^2+f^2}=\frac{W}{R^2}>0.
\end{equation}

Since $f/R=\sin\theta$, the six simple zeros of $f$ in one circuit
of $\partial D$ are exactly the crossings of integer multiples of
$\pi$ by the increasing angle $\theta$. Thus $\theta$ increases by
$6\pi$ in one circuit: the normalized pair $(H,f)/R$ winds three
times around the unit circle. With $u$ taken modulo $2\pi$, its lift
satisfies
\[
 \theta(u+2\pi)=\theta(u)+6\pi.
\]
Consequently
\begin{equation}\label{eq:angular-coordinate}
 v(u)=\theta(u)/3,\qquad
 v'(u)>0,\qquad v(u+2\pi)=v(u)+2\pi
\end{equation}
defines an orientation-preserving angular coordinate on the entire
boundary. As $v$ makes one circuit, the normalized coefficient
$f/R=\sin(3v)$ makes three oscillations. On the outer end of
$L_{\rm in}$ we have
\begin{equation}\label{eq:unnormalized-canonical}
 H\,dr+f\,dy
 =R\bigl(\cos(3v)\,dr+\sin(3v)\,dy\bigr).
\end{equation}

Let $\widetilde D$ be the meridian disk enlarged by the collar
$L_{\rm in}$, and let $D_{\rm can}$ be a copy with boundary coordinate
$v$. Extend this boundary reparameterization to an
orientation-preserving disk diffeomorphism
\[
 \varphi:\widetilde D\longrightarrow D_{\rm can},\qquad
 (r,u)\longmapsto(r,v(u))\quad\text{on the boundary collar}.
\]
Here $r$ is the collar coordinate already in use. The displayed
formula specifies the map near the boundary; choose any smooth disk
extension in the interior. Put
\[
 V_{\rm can}=D_{\rm can}\times S^1_y,\qquad
 \Phi=\varphi\times\mathrm{id}_{S^1_y}:
 V\cup L_{\rm in}\longrightarrow V_{\rm can}.
\]
If $\widetilde\xi_E$ denotes the contact structure on $V\cup L_{\rm in}$,
define
\[
 \eta_E=\Phi_*\widetilde\xi_E,\qquad
 (\eta_E)_{\Phi(p)}=d\Phi_p\bigl((\widetilde\xi_E)_p\bigr).
\]
Thus $\eta_E$ is the same plane field expressed on the new solid torus
by the differential of $\Phi$. A defining form transforms by
$(\Phi^{-1})^*$.

The transformed form on the boundary collar is
\eqref{eq:unnormalized-canonical}, with $R$ now regarded as a function
of $v$. Multiply this form by $1/R$ on that collar and extend the
multiplier to a positive smooth function over the interior. Positive
rescaling preserves the cooriented contact structure. After this
normalization, the chosen \emph{canonical contact form} on the whole
boundary collar is
\begin{equation}\label{eq:canonical}
 \alpha_{\rm can}=\cos(3v)\,dr+\sin(3v)\,dy.
\end{equation}

Let $\rho:D_{\rm can}\longrightarrow D_{\rm can}$ be an
orientation-preserving diffeomorphism whose collar expression is
$(r,v)\mapsto(r,v+2\pi/3)$, and set
$\widehat\rho=\rho\times\mathrm{id}_{S^1_y}$.
This collar rotation preserves the contact form
\eqref{eq:canonical}. Define two further contact structures on the
same marked solid torus by
\[
 \eta_C=(\widehat\rho)_*\eta_E,\qquad
 \eta_B=(\widehat\rho^2)_*\eta_E.
\]
Their defining forms can be taken to equal $\alpha_{\rm can}$ on
the same boundary collar. The rotation moves every endpoint two
positions forward in the cyclic labeling. Applying it once to
$E=(16)(25)(34)$ gives $C$, and applying it twice gives $B$, where
\begin{equation}\label{eq:matchings}
 B=(12)(36)(45),\qquad C=(14)(23)(56).
\end{equation}
This displacement by two marked points preserves all region signs.
Figure~\ref{fig:meridional-matchings} shows the three signed matchings
with the same boundary labels.

\begin{figure}[tbp]
\centering
\begin{tikzpicture}[
  x=1cm,y=1cm,font=\small,
  dividing/.style={draw=black,line width=.85pt,line cap=round},
  meridian path/.style={draw=figureblue,line width=.7pt,
    dash pattern=on 3pt off 1.9pt,->,>=stealth,shorten >=2.4pt},
  region sign/.style={font=\scriptsize,inner sep=0pt},
  endpoint/.style={inner sep=0pt,font=\scriptsize}
]
  \foreach \state/\offset in {E/0,B/4.65,C/9.30}{
    \begin{scope}[xshift=\offset cm]
      \draw[rounded corners=3pt,line width=.45pt] (0,-.86) rectangle (3.30,.86);
      \if E\state
        \foreach \x in {.42,1.65,2.88}
          \draw[dividing] (\x,-.86) -- (\x,.86);
        \foreach \x/\sgn in {.21/-,1.035/+,2.265/-,3.09/+}
          \node[region sign] at (\x,.26) {$\sgn$};
      \else
        \if B\state
          \def\flip{1}\def\origin{0}
        \else
          \def\flip{-1}\def\origin{3.30}
        \fi
        \begin{scope}[shift={(\origin,0)},xscale=\flip]
          \draw[dividing] (.42,.86) .. controls (.42,.23) and (1.65,.23) .. (1.65,.86);
          \draw[dividing] (2.88,.86) .. controls (2.88,.39) and (.42,-.39) .. (.42,-.86);
          \draw[dividing] (1.65,-.86) .. controls (1.65,-.23) and (2.88,-.23) .. (2.88,-.86);
        \end{scope}
        \node[region sign] at (.44,.16) {$-$};
        \node[region sign] at (2.86,-.16) {$+$};
        \if B\state
          \node[region sign] at (1.035,.65) {$+$};
          \node[region sign] at (2.265,-.65) {$-$};
        \else
          \node[region sign] at (2.265,.65) {$-$};
          \node[region sign] at (1.035,-.65) {$+$};
        \fi
      \fi
      \draw[meridian path] (0,0) -- (3.30,0);
      \foreach \x/\topmark/\bottommark in {.42/1/6,1.65/2/5,2.88/3/4}{
        \fill (\x,.86) circle (1.35pt);
        \fill (\x,-.86) circle (1.35pt);
        \node[endpoint,above=2.5pt] at (\x,.86) {$\topmark$};
        \node[endpoint,below=2.5pt] at (\x,-.86) {$\bottommark$};
      }
      \foreach \x/\y in {0/0,1.035/.86,2.265/.86,3.30/0,2.265/-.86,1.035/-.86}
        \draw[fill=white,line width=.45pt] (\x,\y) circle (1.55pt);
      \node[left=2pt,inner sep=1pt] at (0,0) {$P$};
      \node[right=2pt,inner sep=1pt] at (3.30,0) {$Q$};
      \node at (1.65,-1.40) {$\state$};
    \end{scope}
  }
\end{tikzpicture}
\caption{The three meridional matchings, drawn up to disk isotopy in the
common boundary marking.
Solid arcs divide the disk into regions marked by the sign of
$\alpha(\partial_y)$. Open circles mark the six intersections of the
torus sutures with the disk boundary. The dashed arc from $P$ to $Q$
crosses all three dividing arcs in $E$, only $(36)$ in $B$, and only
$(14)$ in $C$. In either alternative its product with $S^1_y$ gives the
replacement annulus.}
\label{fig:meridional-matchings}
\end{figure}
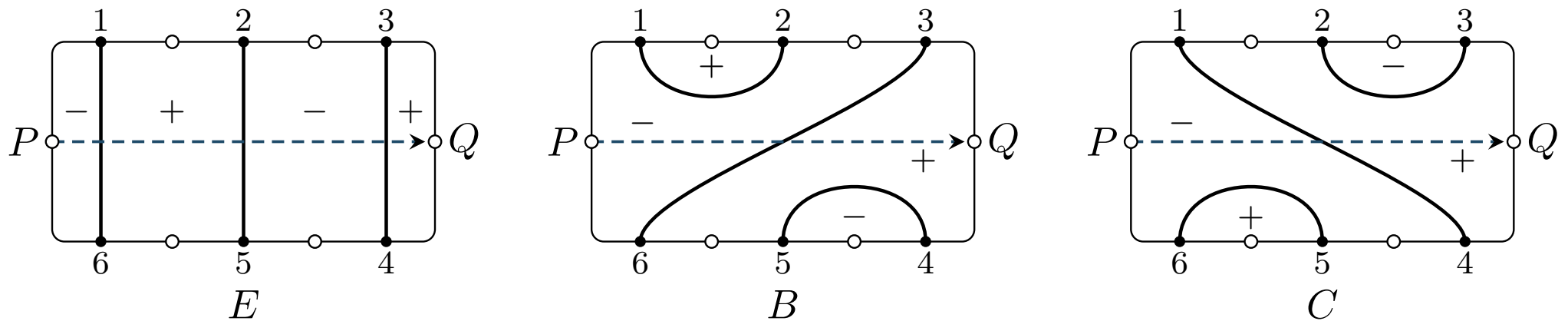

To recover the contact form $\alpha_{\rm old}$ on the outer boundary
collar, attach a second, \emph{restoring collar}
\[
 L_{\rm out}=[0,1]_r\times(S^1_u\times S^1_y).
\]
Use \eqref{eq:collar} with the same $H,f,b$ and with $k(0)=0$,
$k(1)=1$. On its inner end collar the contact form is $H\,dr+f\,dy$;
on its outer end collar it is $\alpha_{\rm old}$.
Attach the inner end to $\partial V_{\rm can}$ using the inverse
boundary reparameterization from Step~2 and restoring the factor $R$
in \eqref{eq:unnormalized-canonical}.
These contact forms then define the same cooriented contact structure
on the full overlap. This attachment is made \emph{identically for
$E$, $B$ and $C$}, since all three fillings have the same defining
contact form $\alpha_{\rm can}$ on their boundary collars.
The outer end of $L_{\rm out}$ is identified with the fixed original
collar of $\partial V$ (translating the radial coordinate if necessary).

\subsection*{Step 3. The relative contact structures and the returning original state}
Choose once and for all a product diffeomorphism
\begin{equation}\label{eq:onej}
 j:(V_{\rm can}\cup L_{\rm out})\longrightarrow V
\end{equation}
which is the prescribed identification of the actual outer collars.
The enlarged meridian disk can be identified with the original meridian
disk by an orientation-preserving diffeomorphism with this prescribed
collar restriction. The $S^1_y$ coordinate is unchanged. In particular
the annular end neighborhoods, meridian, and longitude are fixed in
the marking. Use this same $j$ for all three fillings.

The returning $E$ state has only added radial endpoint extensions to its
three dividing arcs. A disk isotopy relative to the outer collar makes
those signed arcs coincide with the original signed arcs. Incorporate
this one isotopy into $j$ for all three fillings. It does not alter the
boundary matching or the points $P,Q$.

For $K\in\{E,B,C\}$, let $\xi_K$ be the contact structure on $Y$
obtained by attaching $L_{\rm out}$ to $(V_{\rm can},\eta_K)$,
mapping the result to $V$ by $j$, and retaining the original
$\xi$ on the exterior. All three structures agree with $\xi$ on
the fixed outer collar and outside $V$. Thus conclusion
\textup{(\ref{replacement:relative})} holds for both candidate
structures $\xi_B,\xi_C$.

\begin{claim}[The returning original state]\label{clm:returning-E}
The contact structures $\xi_E$ and $\xi$ are isotopic through contact
structures fixed on the outer collar of $V$ and on $Y\setminus V$.
\end{claim}
\begin{proof}
Choose an $S^1_y$-invariant defining form $\widehat\alpha_E$ for
$\xi_E|_V$. On the outer collar it is a positive invariant multiple
of the original form $\alpha_E=\sin x\,dy+\cos x\,dz$.
Extend the inverse of this multiplier positively to $V$ and rescale
$\widehat\alpha_E$, so that $\widehat\alpha_E=\alpha_E$ on that
collar. By the choice of $j$, their $dy$ coefficients have the same
signed regular zero set. Lemma~\ref{lem:fixedf} gives a contact path
between their kernels, fixed on the collar. Extending this path by
the constant contact structure $\xi$ outside $V$ proves the claim.
\end{proof}

This contact isotopy identifies the returning $E$ structure with the
original one. The common product marking remains the single map
$j$ in \eqref{eq:onej} for all three fillings.

\subsection*{Step 4. The six-point relation}
Use $\F$ coefficients for all Floer groups. Let $\Gamma_{\rm can}$ be
the six-component dividing set on $\partial V_{\rm can}$, and let
\[
 e_K^0=\EH(V_{\rm can},\Gamma_{\rm can},\eta_K)
 \in\SFH(-V_{\rm can},-\Gamma_{\rm can}),\qquad K\in\{E,B,C\}.
\]
These elements use the same marked boundary.

In the disk contact TQFT of \cite[Section 7.1]{HKMTQFT}, the matchings
$E,B,C$ are the three central-Euler-grading states for six endpoints.
To identify the constructed fillings with these states, align
the signed zero arcs of two invariant realizations by a disk isotopy
fixed on the boundary collar. Lemma~\ref{lem:fixedf} then gives a
contact path fixed on that collar. Thus the signed matching determines
the marked contact element. The six-point relation
\cite[Lemma 7.4]{HKMTQFT} gives
\begin{equation}\label{eq:HKMrelation}
 e_E^0=e_B^0+e_C^0.
\end{equation}

\subsection*{Step 5. The annulus and contact vector field}
\emph{Conclusion \textup{(\ref{replacement:annulus})}.}
For $K=B,C$, let $\gamma_K$ be the proper arc represented by the dashed
arc from $P$ to $Q$ in Figure~\ref{fig:meridional-matchings}, transported
to the marked meridian disk. Choose it to agree with the original
horizontal arc on the fixed end collars. The arc $\gamma_B$ crosses
only the dividing chord $(36)$, and $\gamma_C$ crosses only $(14)$;
each crossing is transverse. The added radial collars extend the ends
without creating further crossings.

Orient $\gamma_K$ from $P$ to $Q$ and set
$A'_K=\gamma_K\times S^1_y$, with the product orientation, in the
common marking of Step~3. Since $\gamma_K$ is properly embedded
in the marked disk, $A'_K$ is a properly embedded annulus in $V$.
The agreement of $\gamma_K$ with the horizontal arc on both end collars gives
$A'_K=A_0$ on both annular end collars.
Any two proper embedded arcs in a disk with these fixed end collars
are isotopic relative to the collars. Taking the product with $S^1_y$
therefore gives a smooth isotopy $A'_K\simeq A_0$ with the required
relative framing. This proves conclusion
\textup{(\ref{replacement:annulus})} for both candidate annuli.

\emph{Conclusion \textup{(\ref{replacement:positive})}.}
Choose an $S^1_y$-invariant contact form $\widehat\alpha_K$ for the alternative
filling that equals the original model form
$\alpha_E=\sin x\,dy+\cos x\,dz$ on the fixed outer collar.
Such a form exists because the contact structures agree there: multiply
any invariant defining form by the required positive invariant function
on that collar and extend the multiplier positively over the interior.
This normalization will give the exact end vector field required in
conclusion \textup{(\ref{replacement:ends})}.

Fix $K\in\{B,C\}$ and abbreviate $A'_K$ to $A'$.
Parametrize $\gamma_K$ by $s\in[-\varepsilon,\varepsilon]\subset(-\pi/2,\pi/2)$, increasing
from $P$ to $Q$, with its unique dividing-arc crossing at $s=0$.
On sufficiently small end collars choose
\[
 s=x\quad\text{at the left end},\qquad
 s=x-2\pi\quad\text{at the right end}.
\]
These prescriptions extend to a smooth increasing parameter along the arc,
since they take negative values on the left and positive values on the right.
The pullback of $\widehat\alpha_K$ to $A'$ has the form
\[
 \beta=\widehat\alpha_K|_{TA'}=p(s)\,dy+q(s)\,ds,
\]
where $p$ changes from negative to positive at its unique simple zero
$s=0$. Since $\varepsilon<\pi/4$, the function
\[
 r(s)=\tan s
\]
is smooth on the whole parameter interval, satisfies $r(0)=0$ and
$r'(s)>0$, and agrees with $\tan x$ on both end collars.

Set $H|_{A'}=p/\tan s$, using the smooth positive extension at $s=0$
from Lemma~\ref{lem:annulus}. On both end collars,
\[
 H=\frac{\sin x}{\tan s}=\cos x.
\]
Extend $H$ to equal $\cos x$ on the full three-dimensional end
neighborhoods and smoothly elsewhere. Let $X_K$ be its contact
Hamiltonian vector field with respect to $\widehat\alpha_K$.
Thus $X_K$ is a contact vector field for
$\xi_K=\ker\widehat\alpha_K$ near $A'_K$.
Lemma~\ref{lem:annulus} and formula~\eqref{eq:ham} show that $X_K$
is positively transverse
to $A'$. This test depends only on the restriction of $H$ and its
tangential differential along $A'$, so the chosen extension
preserves transversality. Moreover,
\[
 \widehat\alpha_K(X_K)|_{A'_K}=H|_{A'_K}>0.
\]
Thus $X_K$ is also positively transverse to $\xi_K$ along $A'_K$,
and its dividing set on $A'_K$ is empty. This proves conclusion
\textup{(\ref{replacement:positive})}.

\emph{Conclusion \textup{(\ref{replacement:ends})}.}
On each end neighborhood, the contact form is $\alpha_E$ and
$H=\alpha_E(\partial_z)$. Uniqueness of the contact Hamiltonian
vector field therefore gives $X_K=\partial_z$ on that full
neighborhood, proving conclusion \textup{(\ref{replacement:ends})}.

\subsection*{Step 6. One gluing map and nonvanishing}
The exterior of $V$ is connected, since every exterior component meets
the connected boundary torus. Remove a fixed Darboux ball from this
exterior. Attach the common restoring collar $L_{\rm out}$ to
$V_{\rm can}$, identify the resulting solid torus with $V$ by the
common map $j$, and attach the fixed punctured exterior. These
attachments define one HKM gluing map
\[
 \Phi_{\rm tot}:
 \SFH(-V_{\rm can},-\Gamma_{\rm can})\longrightarrow\HFhat(-Y).
\]
The attached complement is connected and meets the puncture boundary,
so there is no isolated component and no additional tensor factor.

For the closed contact structures $\xi_K$ defined in Step~3,
contact naturality
\cite[Theorem 1.1]{HKMTQFT} gives
\[
 \Phi_{\rm tot}(e_K^0)=c(\xi_K).
\]
Claim~\ref{clm:returning-E} identifies $c(\xi_E)$ with $c(\xi)$.
Applying this single linear map to
\eqref{eq:HKMrelation} therefore gives
\begin{equation}\label{eq:closedrelation}
 c(\xi)=c(\xi_B)+c(\xi_C).
\end{equation}
Since $c(\xi)\ne0$, at least one term on the right is nonzero.
Choose $K\in\{B,C\}$ with $c(\xi_K)\ne0$ and set
\[
 \xi'=\xi_K,\qquad A'=A'_K,\qquad X=X_K.
\]
This choice gives $c(\xi')\ne0$, proving
\textup{(\ref{replacement:nonzero})}. By Step~3, $\xi'$ equals
$\xi$ outside $V$ and on the fixed open boundary collar, proving
\textup{(\ref{replacement:relative})} for the selected structure.
Step~5 gives the properly embedded annulus $A'$ with its stated
agreement and isotopy relative to both end collars, proving
\textup{(\ref{replacement:annulus})}. The same step gives the
$\xi'$-contact vector field $X$, positively transverse to $A'$
and to $\xi'$ along $A'$, with empty dividing set, proving
\textup{(\ref{replacement:positive})}. It also gives
$X=\partial_z$ on the full three-dimensional end neighborhoods,
proving \textup{(\ref{replacement:ends})}. All five conclusions
therefore hold for the same triple $(\xi',A',X)$.
\end{proof}

\section{Realizing the rotation model}

We prove Proposition~\ref{prop:P} by constructing a surface one-form satisfying the desired local model and then realizing it through a relative isotopy. 
We first establish a relative realization lemma and then use a positive area-form extension to construct the required one-form on the remaining positive regions. 

\subsection{Relative realization}
\label{sec:realization}

The following lemma realizes a new surface one-form and the associated contact structure on a flow neighborhood. 
The proof uses the Hamiltonian identity established in Lemma~\ref{lem:ham}. 

\begin{lemma}[Relative realization of an invariant contact neighborhood]
\label{lem:realization}
Let $(Y,\xi)$ be a positive cooriented closed contact three-manifold, and let $i_0:F\hookrightarrow(Y,\xi)$ be a compact oriented convex surface with positively transverse boundary.
Choose a positively transverse contact vector field $X_0$ whose dividing set is disjoint from $\partial F$, and use its flow to identify an invariant neighborhood of $i_0(F)$ with
\[
 \xi=\ker(f\,dt+\beta),\qquad f\,d\beta+\beta\wedge df>0,
\]
where $f$ and $\beta$ are defined on $F$, $i_0(F)$ is the zero slice, and $X_0=\partial_t$ in these coordinates. 
Set $\Gamma=f^{-1}(0)$. 
Suppose that a smooth one-form $\beta_*$ on $F$ satisfies
\[
 Q_*=f\,d\beta_*+\beta_*\wedge df > 0
\]
and agrees with $\beta$ on a neighborhood of the union of the closures of all components of $F\setminus\Gamma$ that meet $\partial F$. 
Then there exist an isotopy of embeddings $i_s$, fixed near $\partial F$, and a flow chart
\[
 \Psi:F\times(-\delta,\delta) \longrightarrow Y, \qquad \Psi(p,0)=i_1(p), 
\]
such that the following properties hold:
\begin{enumerate}
\item $\Psi^*\xi=\ker(f\,dt+\beta_*)$. 
\item The vector field $X=\Psi_*\partial_t$ is a contact vector field in the flow neighborhood and is positively transverse to $i_1(F)$. 
Moreover, its dividing set pulls back under $i_1$ exactly to $\Gamma$, and its closed positive and negative regions pull back exactly to $\{f\geq0\}$ and $\{f\leq0\}$, respectively. 
\item The vector field $X$ agrees with $X_0$ on a full ambient neighborhood of $i_0(\partial F)$.  
\end{enumerate}
\end{lemma}
\begin{proof}
The inequality $Q_*>0$ is the dividing-set criterion from convex surface theory.
In the present notation, this can be checked directly: on the auxiliary product $F\times(-\delta,\delta)$, the one-form $\alpha_*=f\,dt+\beta_*$ satisfies
\[
 \alpha_*\wedge d\alpha_*=dt\wedge Q_*>0.
\]
Since $\alpha_*$ is independent of $t$, the vector field $\partial_t$ is a contact vector field positively transverse to
the zero slice, and its dividing set on that slice is $\{f=0\}$. 
Thus the oriented singular foliation defined by $\beta_*$ is divided by $\Gamma$. 
By the hypothesis that $\beta_*=\beta$ on a neighborhood of the closures of all components of $F\setminus\Gamma$ meeting $\partial F$, this foliation agrees with the original one on those components. 
The relative flexibility theorem \cite[Theorem~2.2]{EVHM} then gives an isotopy $i_s$, fixed near $\partial F$, such that the characteristic foliation of $i_1(F)$, pulled back by $i_1$, is the oriented singular foliation defined by $\beta_*$. 

\emph{Comparison of the defining one-forms, including at their zeros.}
Choose a contact form $\alpha$ for $\xi$ that agrees with $f\,dt+\beta$ on the initial invariant neighborhood, and set $\beta_{\rm act}=i_1^*\alpha$. 
Let $Z$ denote the common zero set of $\beta_{\rm act}$ and $\beta_*$. 
On $F\setminus Z$, the equality of their oriented kernels determines a unique smooth positive function $h_0$ such that $\beta_{\rm act}=h_0\beta_*$. 
We claim that $h_0$ extends to a smooth positive function $h$ on all of $F$. 

Fix a point $p\in Z$. 
Since $Q_*(p)=f(p)d\beta_*(p)>0$, we have $d\beta_*(p)\ne0$. 
Moreover, $\beta_{\rm act}(p)=0$ implies that
\[
 di_1(T_pF)=\xi_{i_1(p)}.
\]
This means that $i_1(F)$ is tangent to $\xi$ at the single point $i_1(p)$; it does not mean that $i_1(F)$ is tangent to $\xi$ on an open set. 
Since $d\alpha|_{\xi}$ is nondegenerate, we obtain
\[
 d\beta_{\rm act}(p)=(i_1^*d\alpha)(p)\ne0.
\]
Thus, there is no conflict with the nonintegrability of the contact plane field. 

Write locally
\[
 \beta_*=a\,dx+b\,dy,\qquad
 \beta_{\rm act}=A\,dx+B\,dy.
\]
After interchanging $x$ and $y$ if necessary, we may assume that $da(p)\ne0$. 
Shrinking the neighborhood if necessary, we may assume that $a$ is a submersion and that $\{a=0\}$ is a coordinate hypersurface. 
The function $A$ vanishes on $\{a=0\}$: away from $Z$, this follows from $A=h_0a$, while on $Z$, both one-forms vanish.
By Hadamard's lemma, there exists a smooth function $h_1$ such that $A=ah_1$. 
On the dense set $\{a\ne0\}$, we have $h_1=A/a=h_0$, and hence $B=h_1b$ there. 
By continuity, the latter identity holds throughout the neighborhood, so that $\beta_{\rm act}=h_1\beta_*$ everywhere. 
At each point $p\in Z$, differentiating $\beta_{\rm act}=h_1\beta_*$ and using $\beta_*(p)=0$ gives 
\[
 d\beta_{\rm act}(p)=h_1(p)d\beta_*(p).
\]
Since both two-forms are nonzero at $p$, we obtain
$h_1(p)\ne0$. 
Since $h_1$ agrees with the positive function $h_0$ on
$\{a\ne0\}$, we conclude that $h_1>0$ throughout the neighborhood. 
On the rest of $F\setminus Z$, uniqueness of the proportionality factor implies that $h_1=h_0$ as well. 
By uniqueness on $F\setminus Z$ and continuity, these local extensions agree on overlaps and hence define a global smooth positive function $h$ satisfying 
\begin{equation}\label{eq:realization-multiplier}
 \beta_{\rm act}=h\beta_*.
\end{equation}
Near $\partial F$, the embedding and the one-forms are unchanged, so $h=1$ there. 

\emph{Choosing a transverse contact vector field with the required boundary behavior.}
Define the restriction of the Hamiltonian to the realized surface by
\[
 H\circ i_1=hf.
\]
Let $H_0=\alpha(X_0)$ on a neighborhood of $i_0(\partial F)$. 
In the initial invariant chart, we have $H_0=f$. 
Since $i_1=i_0$ and $h=1$ near $\partial F$, the prescribed restriction $hf$ agrees with $H_0\circ i_1$ there. 
Consequently, $hf$ extends to a smooth function $H$ on a neighborhood of $i_1(F)$ that agrees with $H_0$ on a full ambient neighborhood of $i_0(\partial F)$. 
To see this, one can start with an arbitrary smooth extension and interpolate it with $H_0$ using a cutoff function supported in an ambient neighborhood of the boundary where their restrictions to $i_1(F)$ agree. 
This interpolation preserves the identity $H\circ i_1=hf$ everywhere. 

The contact Hamiltonian vector field $X_H$ is positively transverse to $i_1(F)$, since Lemma~\ref{lem:ham} gives
\[
 i_1^*\bigl(\iota_{X_H}(\alpha\wedge d\alpha)\bigr)
 =(hf)d(h\beta_*)+h\beta_*\wedge d(hf)
 =h^2Q_*>0.
\]
In particular, this calculation depends only on $H\circ i_1$ and its differential along $F$; the choice of extension away from $i_1(F)$ does not affect it. 
On the full ambient neighborhood of $i_0(\partial F)$ where $H=H_0$, uniqueness of the contact Hamiltonian vector field implies that $X_H=X_0$. 
Moreover, $\alpha(X_H)\circ i_1=hf$ and $h>0$, so
\begin{equation}\label{eq:realization-signed-sets}
 \begin{aligned}
 i_1^{-1}\bigl(\{q\in i_1(F):\alpha(X_H)(q)=0\}\bigr)
 &=\{f=0\},\\
 i_1^{-1}\bigl(\{q\in i_1(F):\pm\alpha(X_H)(q)\geq0\}\bigr)
 &=\{\pm f\geq0\}.
 \end{aligned}
\end{equation}
These are precisely the original signed subsets, since $\alpha(X_0)=f$ in the initial invariant chart.
Thus, the equalities hold as equalities of subsets of the fixed parameter surface $F$, not merely up to isotopy; the corresponding embedded subsets of $Y$ move with the surface. 

\emph{Identifying the contact structure in flow coordinates.}
Let $\phi_t$ denote the local flow of $X_H$.
Since $X_H$ is a contact vector field, we have $\phi_t^*\alpha=a_t\alpha$ for a positive function $a_t$.
By transversality and compactness of $F$, for sufficiently small $\delta>0$ the map
\[
 \Psi(p,t)=\phi_t(i_1(p))
\]
defines a flow chart.
To determine the contact structure throughout this flow neighborhood, we compute the pullback of $\alpha$ on vectors tangent to $F$ and on the flow direction. 
For $v\in T_pF$, we have
\[
 (\Psi^*\alpha)_{(p,t)}(v)
 =a_t(i_1(p))\,\beta_{\rm act}(v)
 =a_t(i_1(p))h(p)\,\beta_*(v).
\]
For the flow direction, we use $d\phi_t(X_H)=X_H\circ\phi_t$, which expresses the invariance of a vector field under its own flow. 
Thus,
\[
 \begin{aligned}
 (\Psi^*\alpha)_{(p,t)}(\partial_t)
 &=\alpha_{\phi_t(i_1(p))}
   \bigl(X_H(\phi_t(i_1(p)))\bigr)\\
 &=(\phi_t^*\alpha)_{i_1(p)}
   \bigl(X_H(i_1(p))\bigr)\\
 &=a_t(i_1(p))h(p)f(p).
 \end{aligned}
\]
Combining these two evaluations, we obtain
\begin{equation}\label{eq:realization-flow-pullback}
 \Psi^*\alpha=a_t(i_1(p))h(p)(f\,dt+\beta_*).
\end{equation}
Since no division by $f$ is involved, this identity also holds along $\Gamma$. 
The factor $a_t(i_1(p))h(p)$ is positive, and hence
\[
 \Psi^*\xi=\ker(f\,dt+\beta_*).
\]
Moreover, the defining form on the right-hand side is independent of $t$, so $\partial_t$ is a contact vector field in these coordinates. 
Taking $X=X_H=\Psi_*\partial_t$, the transversality, the exact signed-set identities, and the boundary agreement established above give all three conclusions. 
\end{proof}

\subsection{Extending the boundary one-form}\label{sec:area-extension}

\begin{lemma}[Relative positive area-form extension]
\label{lem:area}
Let $Q$ be a compact connected oriented surface with nonempty boundary. 
Suppose that a one-form $\lambda_\partial$ is defined on a neighborhood of $\partial Q$ and satisfies $d\lambda_\partial>0$. 
If
\[
 B=\int_{\partial Q}\lambda_\partial>0,
\]
then there exists a one-form $\lambda$ on $Q$ that agrees with $\lambda_\partial$ on a (possibly smaller) boundary collar and satisfies $d\lambda>0$ everywhere. 
\end{lemma}
\begin{proof}
Choose a boundary collar $C$ sufficiently small so that 
\[
\int_C d\lambda_\partial<B. 
\]
Extend $\lambda_\partial$ from a smaller collar to a one-form $\lambda_0$ on $Q$. 
Choose a positive area form $\omega$ on $Q$ that agrees with $d\lambda_0$ on a still smaller collar and has total integral $B$. 
Such a form can be constructed by first choosing a sufficiently small positive area form on the remaining interior and then adding a nonnegative bump form supported in the interior to adjust its total integral. 

The two-form $\omega-d\lambda_0$ has compact support in $\Int Q$, and its integral vanishes by Stokes' theorem. 
Since integration induces an isomorphism $H_c^2(\Int Q;\R)\cong\R$, there exists a compactly supported one-form $\eta$ on $\Int Q$ such that $\omega-d\lambda_0=d\eta$. 
Then $\lambda=\lambda_0+\eta$ has the required properties. 
\end{proof}

\subsection{Proof of relative preparation}\label{sec:preparation-proof}

\begin{proof}[Proof of Proposition~\ref{prop:P}]
We first construct a one-form $\beta_*$ on the fixed parameter surface $F$, while keeping the function $f$ fixed. 
We then use Lemma~\ref{lem:realization} to realize this new one-form by an isotopy of the surface within the original contact structure. 

\paragraph{Step 1. Fixing $f$ and defining the new form near $A$}
The flow of $X_0$ provides an invariant neighborhood of $i_0(F)$ in which 
\[
 \alpha=f\,dt+\beta,\qquad
 f\,d\beta+\beta\wedge df>0.
\]
Here, $\alpha$ denotes the ambient contact form, whereas $\beta=i_0^*\alpha$ is a one-form on $F$. 
The initial dividing set and signed regions on $F$ are given by
\[
 \Gamma_0=\{f=0\},\qquad
 R_{+,0}=\{f\geq0\},\qquad R_{-,0}=\{f\leq0\}.
\]
We will replace $\beta$ with $\beta_*$ while keeping $f$ fixed. 

Choose a slightly larger annulus $\widehat U$ such that $U\subset\Int\widehat U$, with $\widehat U$ still disjoint from $N$ and from all unselected negative annuli. 
Choose coordinates $(x,y)$ on an open annular neighborhood of $\widehat U$, where $x$ is an interval coordinate transverse to the annular core circles and $y\in S^1$ is the angular coordinate. 
Thus, $x$ is a radial coordinate on the annulus, not an angular variable. 
Choose the coordinates so that $dx\wedge dy$ agrees with the positive surface orientation and 
\[
 U=[-\varepsilon,2\pi+\varepsilon]\times S^1_y,\qquad
 0<\varepsilon<\widehat\varepsilon<\pi/4,\qquad
 \widehat U=[-\widehat\varepsilon,2\pi+\widehat\varepsilon]
 \times S^1_y, 
\]
with the two components of $\partial A$ lying at $x=\pi/2$ and $x=3\pi/2$. 
The region between these circles is negative, while the two end regions of $U$ are positive. 
Zero is a regular value of $f$, since the contact inequality implies that $\beta\wedge df>0$ along $f^{-1}(0)$. 
By the annular configuration, the only components of $\Gamma_0$ in $\widehat U$ are the two boundary circles of $A$, and $f$ is negative between them and positive on the two end regions. 
Thus $f$ and $\cos x$ have the same zero set and the same signs on its complement.
Consequently, the quotient $f/\cos x$ extends smoothly and positively across the two zero curves. 
Denote this extension by $g$, so that $f=g\cos x$, and define
\begin{equation}\label{eq:target}
 \beta_*=g\sin x\,dy
\end{equation}
on a neighborhood of $\widehat U$. 
The terms involving $dg$ cancel, giving 
\begin{equation}\label{eq:targetpositive}
f\,d\beta_*+\beta_*\wedge df=g^2\,dx\wedge dy>0. 
\end{equation}
Outside this annular neighborhood, set $\beta_*=\beta$ near $N$ and near every unselected negative annulus. 
We also set $\beta_*=\beta$ on each positive component disjoint from $\widehat U$.
It remains to define $\beta_*$ on the positive components adjacent to $\widehat U$ by interpolating between the new model form and the original form. 

\paragraph{Step 2. Boundary data and extension over the positive regions}
Let $P$ be a component of $\overline{R_{+,0}\setminus\widehat U}$ whose boundary meets $\partial\widehat U$.
The boundary of $P$ consists of some end circles of $\widehat U$ and some unselected dividing curves.
Passing from $U$ to $\widehat U$ only removes additional boundary collars from the adjacent positive regions.
In particular, it does not change their connectivity or eliminate any unselected dividing boundary.
Consequently, the outlet condition guarantees that $\partial P$ contains at least one unselected dividing curve. 

For each unselected dividing boundary component of $P$, remove a sufficiently small collar from its positive side.
More explicitly, use $s=f\geq0$ as the inward collar coordinate and remove the region $0\leq s<\delta$.
Choose these collars to be pairwise disjoint and disjoint from the model-end collars. 
Denote the remaining compact connected surface by $Q$. 
Its boundary consists of the end circles inherited from $\partial\widehat U$ and the truncation circles $s=\delta$. 
We define the boundary one-form $\lambda_\partial$ for Lemma~\ref{lem:area} separately on these two types of boundary. 
Near an end circle of $\widehat U$, define
\[
 \lambda_\partial=\frac{g\sin x\,dy}{g\cos x}
 =\tan x\,dy.
\]
This form is defined on a full collar in $Q$, using the annular coordinates from Step~1. 
Near a truncation circle, define instead
\[
 \lambda_\partial=\frac{\beta}{f},
\]
where $\beta$ is the original one-form on $F$. 
We reserve the notation $\lambda$ for the extension over $Q$ constructed below. 
Both types of boundary data have positive exterior derivative. 
For the first, we have
\[
 d(\tan x\,dy)=\sec^2x\,dx\wedge dy>0;
\]
for the second, positivity follows from
\begin{equation}\label{eq:liouvillepositive}
 d\left(\frac{\beta}{f}\right) = \frac{f\,d\beta+\beta\wedge df}{f^2}>0. 
\end{equation}

The remaining hypothesis of Lemma~\ref{lem:area} is that 
\[
 \int_{\partial Q}\lambda_\partial>0.
\]
This is where the outlets play a role. 
On the end circle $x=-\widehat\varepsilon$, the boundary orientation induced by $Q$ is the positive $y$ direction, whereas on $x=2\pi+\widehat\varepsilon$, it is the negative $y$ direction. 
Consequently, each end circle contributes the same fixed negative amount, 
\[
 -\tan\widehat\varepsilon\int_{S^1}dy<0. 
\]
These negative contributions must be compensated for by the positive contributions from the truncation circles, where the original form $\beta/f$ is used. 

To compute this contribution, choose an angular coordinate $\vartheta$ along the chosen unselected dividing curve, increasing in its boundary orientation. 
Since $s=f$ is an inward-pointing coordinate, the positive surface orientation near this curve is given by $-ds\wedge d\vartheta$. 
Write
\[
 \beta=a(s,\vartheta)\,d\vartheta+b(s,\vartheta)\,ds. 
\]
At $s=0$, the original contact inequality reduces to 
\[
 \beta\wedge ds
 =-a(0,\vartheta)\,ds\wedge d\vartheta>0,
\]
so $a(0,\vartheta)>0$ at every point of the dividing curve. 
The boundary orientation of $Q$ along $s=\delta$ is also the positive $\vartheta$ direction. 
Hence
\begin{equation}\label{eq:outletintegral}
 \int_{s=\delta}\lambda_\partial
 =\frac1\delta\int_{S^1}a(\delta,\vartheta)\,d\vartheta
 \longrightarrow+\infty
 \qquad(\delta\longrightarrow0^+). 
\end{equation}
Indeed, the numerator converges to $\int_{S^1}a(0,\vartheta)\,d\vartheta>0$. 
Since $P$ has at least one outlet, we can choose $\delta$ sufficiently small so that the sum of the contributions from the truncation circles exceeds the absolute value of the contributions from the model-end circles. 
Thus $\int_{\partial Q}\lambda_\partial>0$. 

Apply Lemma~\ref{lem:area} to each of these compact surfaces $Q$. 
We obtain a one-form $\lambda$ on $Q$ that agrees with $\lambda_\partial$ on a possibly smaller boundary collar and satisfies $d\lambda>0$. 
Define
\[
 \beta_*=f\lambda\quad\text{on }Q. 
\]
By the prescribed boundary values of $\lambda_\partial$, this agrees with $g\sin x\,dy$ near the model-end boundaries and with $\beta$ near the truncation circles. 
Set $\beta_*=\beta$ on the strips between the unselected dividing curves and their truncation circles, as well as on every component of $F$ not involved in this construction.
Since these definitions agree on full collars of the corresponding boundary components, they fit together smoothly. 
On each $Q$, we have
\[
 f\,d\beta_*+\beta_*\wedge df=f^2d\lambda>0. 
\]
Together with \eqref{eq:targetpositive} and the fact that
$\beta_*=\beta$ on the remaining pieces, this proves
\begin{equation}\label{eq:globaltarget}
 f\,d\beta_*+\beta_*\wedge df>0\quad\text{on }F. 
\end{equation}
In particular, $\beta_*=\beta$ on a neighborhood of the closed protected collars $N$. 

\paragraph{Step 3. Realization and conclusions
\textup{(\ref{prep:boundary}) and (\ref{prep:regions})}}
The components of $F\setminus\{f=0\}$ that meet $\partial F$ are precisely the interiors, relative to $F$, of the protected negative collars.
Their closures are the $N_i$, and $\beta_*=\beta$ on a neighborhood of their union. 
Thus, \eqref{eq:globaltarget} and this equality verify the
hypotheses of Lemma~\ref{lem:realization}.
The lemma gives an isotopy $i_s$, fixed near $\partial F$,
and a flow chart $\Psi$ near $i_1(F)$ such that
\[
 \Psi^*\xi=\ker(f\,dt+\beta_*),\qquad
 \Psi(p,0)=i_1(p).
\]
Let $X_1$ be the contact vector field corresponding to $\partial_t$ in these flow coordinates. 
It is positively transverse to $i_1(F)$, and the defining one-form $f\,dt+\beta_*$ evaluates on $X_1$ as $f$. 
Consequently, the dividing set of $X_1$ pulls back under $i_1$ exactly to $\{f=0\}$, and its closed positive and negative regions pull back exactly to $\{f\geq0\}$ and
$\{f\leq0\}$, respectively. 
Since $f$ has remained fixed throughout the construction, these are precisely the original subsets of $F$.
This proves \textup{(\ref{prep:regions})} as literal equalities of subsets of $F$. 
This preservation follows from the realization lemma together with the fact that $f$ is unchanged; the role of Lemma~\ref{lem:area} was to construct the extensions with positive exterior derivative in Step~2. 

\paragraph{Step 4. The model and conclusions \textup{(\ref{prep:model}) and (\ref{prep:field})}}
On $\widehat U$, the formulas $f=g\cos x$ and $\beta_*=g\sin x\,dy$ give
\[
 \ker(f\,dt+\beta_*)
 =\ker\bigl(g(\cos x\,dt+\sin x\,dy)\bigr)
 =\ker(\cos x\,dt+\sin x\,dy),
\]
since $g>0$. 
Setting $z=t$, we obtain the rotation model on a sufficiently small flow neighborhood of $i_1(U)$, with $X_1=\partial_t=\partial_z$. 
This proves \textup{(\ref{prep:model})} and the model part of \textup{(\ref{prep:field})}. 
By property~(3) of Lemma~\ref{lem:realization}, $X_1=X_0$ on a neighborhood of $T$, which proves the remaining part of \textup{(\ref{prep:field})}. 

Choose the replacement solid torus inside this flow neighborhood, taking it sufficiently thin so that it intersects the realized surface only along the selected annular patch. 
It is then disjoint from the images of all unselected annuli and all other components of $F$, as required for the subsequent application of Proposition~\ref{prop:L}. 
\end{proof}

\section{Further questions}

As discussed in Remark~\ref{rem:transverse-boundary},
Theorem~\ref{main:surface} does not establish preservation of
the transverse isotopy class of the boundary.
The following question, formulated by Ito and Kawamuro
\cite[Conjecture~2]{IK}, remains beyond the scope of our proof.

\begin{question}[Ito--Kawamuro]
Let $T$ be a positively transverse link in $(S^3,\xi_{\std})$
satisfying
\[
\slink(T)=-\chi_{\max}(T).
\]
Is $T$ transversely isotopic to the closure of a strongly
quasipositive braid?
\end{question}

It is also natural to ask whether the characterization of
strong quasipositivity in terms of Bennequin sharpness extends
to general tight contact three-manifolds.
See \cite[Conjecture~4.1]{Hayden} and
\cite[Problem~1.69, Remark~(5)]{K3}.

\begin{question}
Let $(Y,\xi)$ be a closed tight contact three-manifold,
and let $T$ be a positively transverse link bounding
a Seifert surface $F$ such that
\[
\slink_F(T)=-\chi(F).
\]
Must $T$ be strongly quasipositive, in the sense that
it bounds a Legendrian ribbon in $(Y,\xi)$?
\end{question}
The authors have obtained a negative answer to this question
by constructing a Bennequin-sharp knot in a tight contact
three-manifold that is not strongly quasipositive.
The details will appear in a forthcoming paper.
\bibliographystyle{amsalpha}
\bibliography{bibliography}
\end{document}